\documentclass[12pt,twoside,a4paper]{article}
\usepackage{graphicx}
\usepackage{amsfonts}
\usepackage{bbm}
\usepackage[leqno]{amsmath}
\usepackage{mathrsfs}
\usepackage{fancyhdr}
 \usepackage{mathrsfs,amsmath,amssymb,amsthm}
\usepackage[dvips]{color}
 \usepackage{amssymb}
 \usepackage{amsbsy}
 \usepackage[dvips]{color}
\usepackage{titlesec}
\titlespacing*{\section}{0pt}{0.5\baselineskip}{0.5\baselineskip}
\titleformat*{\subsubsection}{\it}
\titlespacing*{\subsection}{0pt}{0.3\baselineskip}{0.3\baselineskip}
\titlespacing*{\subsubsection}{0pt}{0.3\baselineskip}{0.3\baselineskip}

\allowdisplaybreaks

 \theoremstyle{definition}
\theoremstyle{remark}  \newtheorem{remark}{\noindent\mbox{Remark}}
 \theoremstyle{plain}
 \theoremstyle{plain}\newtheorem{lemma}{\noindent\mbox{Lemma}}
\theoremstyle{plain} \newtheorem{theorem}{\noindent\mbox{Theorem}}
 \theoremstyle{plain}\newtheorem{proposition}{\noindent\mbox{Proposition}}
 \theoremstyle{plain}
\theoremstyle{definition} 

 \def\proof{\noindent{\it Proof.~~}}
 \def\qed{\hfill$\Box$\medskip}
 \def\rto{\rightarrow\infty}
\def\z{\left}
\def\y{\right}
 \def\no{\nonumber}
 \def\mb{\mathbf}

\begin{document}

 \title{{Two-type linear-fractional branching processes in varying environments with asymptotically constant mean matrices}}                

\author{Hua-Ming \uppercase{Wang}\footnote{Email:hmking@ahnu.edu.cn; School of Mathematics and Statistics, Anhui Normal University, Wuhu, 241003, China }~ $\ \&\ $
Huizi \uppercase{Yao}\footnote{Email:yaohuiziyao@163.com; School of Mathematics and Statistics, Anhui Normal University, Wuhu, 241003, China}
}
\date{}
\maketitle%

\vspace{-1cm}

\begin{center}
\begin{minipage}[c]{12cm}
\begin{center}\textbf{Abstract}\quad \end{center}
Consider two-type linear-fractional branching processes in varying environments with asymptotically constant mean matrices.  Let $\nu$ be the extinction time. Under certain conditions, we show that both $P(\nu=n)$ and $P(\nu>n)$ are asymptotically the same as some functions of the products of spectral radii of the mean matrices.
We also give an example for which $P(\nu=n)$ decays with various speeds such as $\frac{c}{n(\log n)^2},$ $\frac{c}{n^\beta},\beta >1$ et al. which are very different from the ones of homogeneous  multitype Galton-Watson processes.

\vspace{0.2cm}

\vspace{0.2cm}

\end{minipage}
\end{center}

\section{Introduction}

Compared with Galton-Watson processes, many new phenomena arise when considering the branching processes in varying environments (BPVE hereafter). For example, for single-type case, Lindvall \cite{lin}
showed that the population size $Z_n$ converges almost surely to a random variable $Z_\infty,$ which may take a positive value with positive probability, that is, in the words of Kersting \cite{ker}, the process may ``fall asleep" at some positive state.
Fujimagari \cite{fuj} showed that the tail probability of the extinction time may behave asymptotically like $(\log n)^{-1}$ or $n^{-\beta}(\beta<1).$  Later, in \cite{ms} Macphee and Schuh discovered that a BPVE may diverge at different exponential rates.
Over the past decades, single-type BPVEs have been extensively studied.   Criteria for almost sure extinction, asymptotics of survival probability,  the
distribution of the population size, Yaglom type limit theorem and many others are known in full generality.
 For details, we refer the reader to  Bhattacharya and Perlman \cite{bp}, Jagers \cite{j}, Kersting \cite{ker}   and Chap.1 of  Kersting and Vatutin \cite{va} and references therein.

There are only a few results about the multitype BPVE and the situation  is less satisfying.
Jones \cite{jon} gave some second moment conditions under which population size of multitype BPVE, normalized by its mean, converges almost surely or in $L^2$ to a random limit. Further studies on $L^p$ and almost sure convergence and the continuity of the limit distribution can be found in Biggins et al. \cite{bcn} and Cohn and Wang \cite{cw}. Recently,  under very general assumptions (arbitrary number of types and general offspring distributions), criteria were provided in  Dolgopyat et al. \cite{dhkp} for
 almost sure extinction and the convergence of the population size (conditioned on survival and normalized by its mean) to an exponential random variable. We notice also that in \cite{dhkp}, the authors  gave in addition the asymptotics of the survival probabilities. Their proof relies on a generalization of the Perron-Frobenius theorem suitable for studying the product of nonhomogeneous nonnegative matrices. Let $\nu$ be the extinction time and $M_k$ be the mean matrix of offspring distribution of individuals of the $(k-1)$-th generation. They showed that \begin{align}\label{sv}
  \frac{1}{C}\Big(\sum_{k=1}^n\lambda_1^{-1}\cdots\lambda_{n}^{-1}\Big)^{-1}<P(\nu>n)<C\Big(\sum_{k=1}^n\lambda_1^{-1}\cdots\lambda_{n}^{-1}\Big)^{-1}, \end{align} where $0<C<\infty$ is certain constant and for $k\ge1,$ $\lambda_k$ is a number associated to the generalized Perron-Frobenius theorem and depends on the mean matrices $M_n,n\ge k$ of the branching process, see \cite[Proposition 2.1 and Lemma 2.2]{dhkp}.
 Inspired by \eqref{sv}, we may consider further the following questions.
 \vspace{-0.3cm}

 \begin{itemize}
      \item[(i)] The number $\lambda_k$ depends on the  mean matrices $M_n,n\ge k$ so that it is hard to compute explicitly. Therefore, instead of $\lambda_k,$ we intend to use $\varrho(M_k),$ the spectral radius of $M_k,$ which is directly computable.
          \vspace{-0.2cm}

 \item[(ii)] Give conditions to ensure a precise asymptotical equivalence  of $P(\nu>n),$ i.e., show that $P(\nu>n)\sim c\z(\sum_{k=1}^{n+1} \varrho(M_1)^{-1}\cdots \varrho(M_{k-1})^{-1}\y)^{-1} \text{ as }n\rto, $ for some constant $0<c<\infty.$
\vspace{-0.2cm}

 \item[(iii)] Consider further the asymptotics of $P(\nu=n).$ We emphasize that even if one knows the asymptotics of $P(\nu>n),$ it is not an easy task to deal with the asymptotics of $P(\nu=n).$
 \end{itemize}
 \vspace{-0.3cm}

Currently, for multitype BPVE with general offspring distribution, we have no idea how to solve the above questions. In this paper, we consider  two-type linear-fractional BPVEs with asymptotically constant mean matrices.
For the linear-fractional setting, the distribution of the extinction time can be written in terms of sum of product $\prod_{k=1}^n M_k$ of the mean matrices.
But $\prod_{k=1}^n M_k$ is hard to be estimated. So we construct some new matrix $A_k$ which may depend on $M_k$ and $M_{k+1}.$ The product $\prod_{k=1}^nA_k$ is related to the approximants of a continued fraction. Thus by some delicate analysis of the continued fractions and the product $\prod_{k=1}^n A_k,$  we can express the asymptotics of extinction time distribution in terms of the product of the spectral radii of $M_k,$ which can be computed explicitly.

We also construct examples for which the mass of the extinction time distribution at $n$ decays with various speeds, for example, $\frac{c}{n(\log n)^2},$ $\frac{c}{n^\beta},\beta>1$ et al. This observation complements  Fujimagari \cite{fuj}, which studied the single-type counterpart.

\noindent{{\bf Outline of the paper.} In Section \ref{mr}, we introduce the  model and conditions, state the main result and give the main body of its proof.
 Then in Section \ref{eg}, we construct a two-type linear-fractional BPVE whose extinction time exhibits asymptotics very different from those of homogeneous Galton-Watson process. Sections \ref{pbp} and \ref{pmp} are devoted to proving the results required for proving the main theorem.

\section{Model and main results} \label{mr}
Suppose that $M_k,k\ge1$ is a sequence of  nonnegative 2-by-2 matrices and $\gamma_k=(\gamma_k^{(1)},\gamma_k^{(2)}), k\ge 1$ is a sequence of nonnegative 2-dimensional row vectors. To avoid the degenerate case, we require that $\forall k\ge1,$ all elements of $M_kM_{k+1}$ are strictly positive and $\gamma_k\ne \mb 0.$
For $\mathbf s=(s_1,s_2)^t\in [0,1]^2$ and $k\ge1,$ let
\begin{align*}
 \mathbf{f}_{k}(\mathbf{s})=(f_{k}^{(1)}(\mathbf{s}),f_{k}^{(2)}(\mathbf{s}))^{t}=\mb1-\frac{M_{k}(\mb1-\mathbf{s})}{1+\gamma_{k}(\mb1-\mathbf{s})}
\end{align*}
which is known as the probability generating function of a linear-fractional distribution. Here and in what follows, $\mb v^t$ denotes the transpose of a vector $\mb v$ and $\mb 1=(\mb e_1+\mb e_2)^t=(1,1)^t,$ with $\mb e_1=(1,0),\mb e_2=(0,1).$

Suppose that $Z_n=(Z_{n,1},Z_{n,2}),n\ge0$ is a stochastic process such that
\begin{align*}
  E\z(\mathbf s^{Z_n}\big|Z_0,...,Z_{n-1}\y)=\z[\mathbf f_{n}(\mb s)\y]^{Z_{n-1}}, n\ge1,
\end{align*}
where $\z[\mathbf f_{n}(s)\y]^{Z_{n-1}}:=\z[f_n^{(1)}(\mb s)\y]^{Z_{n,1}}\z[f_n^{(2)}(\mb s)\y]^{Z_{n,2}}.$
We call the process $Z_n,n\ge0$ a two-type linear-fractional branching process in a varying environment. Matrices $M_k,k\ge1$ are usually referred to as the mean matrices of the branching process.
Denote by $$\nu=\min\{n: Z_{n}=\mathbf 0\}$$ the extinction time of $\{Z_n\}$ which we concern.

 Throughout, we assume $b_k,d_k>0,$  $a_k,\theta_k\ge0,$ $a_k+\theta_k>0,$ $k\ge1$  and put
\begin{align}\label{mg}\forall k\ge1, M_k:=\left( \begin{array}{cc}
  a_k & b_k \\
  d_k &\theta_k \\
 \end{array}\right),\gamma_k:=\mb e_1M_k.\end{align}

 We introduce the following conditions on the number $a_k,b_k,d_k,\theta_k,k\ge1.$

\noindent{\bf(B1)} Suppose that $b,d>0,$  $ a,\theta\ge0$ are some numbers  such that $a+\theta>0,$
$ a_k\rightarrow a,   b_k\rightarrow b ,  d_k\rightarrow d,\theta_k\rightarrow\theta$ as $k\rto$
 and assume further that
\begin{align}
  \sum_{k=2}^\infty|a_k-a_{k-1}|+| b_k- b_{k-1}|+| d_k- d_{k-1}|+|\theta_{k}-\theta_{k-1}|<\infty.\label{abs}
\end{align}

Suppose now condition (B1) holds and for $k\ge1,$ set \begin{align}\label{dta} A_k:=\left( \begin{array}{cc}
  \tilde a_k & \tilde b_k \\
  \tilde d_k &0 \\
 \end{array}\right) \text{ with }
\tilde a_k=a_k+\frac{b_k\theta_{k+1}}{ b_{k+1}}, \tilde b_k= b_k,\tilde  d_k=d_k-\frac{a_k\theta_k}{b_k}.\end{align}
Letting $\Lambda_k=\left(
                 \begin{array}{cc}
                   1 & 0 \\
                   \theta_k/b_k & 1 \\
                 \end{array}
               \right),k\ge1,$ then for $n\ge k\ge1,$ we have
\begin{align}\label{am}
  A_k=\Lambda_k^{-1}M_k\Lambda_{k+1}\text{ and }\mb e_1\prod_{i=k}^n M_i\mb 1=\mb e_1 \prod_{i=k}^n A_i(1,1-\theta_{n+1}/b_{n+1})^t.
\end{align}
\begin{remark}
Note that  the distribution of $\nu$ is formulated by $M_k,$ see \eqref{pngn} and \eqref{pnen} below.
However, the elements of $\prod_{i=k}^n M_i$ are hard to compute and evaluate whereas those of $\prod_{i=k}^n A_i$ are workable because they have some correspondence with continued fractions due to the special structure of the matrices $A_i,i\ge1.$ Therefore, instead of $M_k,$ we will work with $A_k$ below.
\end{remark}
We need in addition the following conditions which are mutually exclusive.

\noindent{\bf(B2)$_{a}$} $\exists k_0>0,$ such that $\frac{\tilde a_k}{\tilde b_k}=\frac{\tilde a_{k+1}}{\tilde b_{k+1}},
\ \frac{\tilde d_k}{\tilde b_k}\ne\frac{\tilde d_{k+1}}{\tilde b_{k+1}},\ \forall k\ge k_0$ and
$$\lim_{k\rto}\frac{\tilde d_{k+2}/\tilde b_{k+2}-\tilde d_{k+1}/\tilde b_{k+1}}{\tilde d_{k+1}/\tilde b_{k+1}-\tilde d_{k}/\tilde b_{k}}\text{ exists.}$$

 \noindent{\bf(B2)$_{b}$} $\exists k_0>0,$ such that $\frac{\tilde a_k}{\tilde b_k}\ne\frac{\tilde a_{k+1}}{\tilde b_{k+1}},
\ \frac{\tilde d_k}{\tilde b_k}=\frac{\tilde d_{k+1}}{\tilde b_{k+1}},\ \forall k\ge k_0$ and
$$\lim_{k\rto}\frac{\tilde a_{k+2}/\tilde b_{k+2}-\tilde a_{k+1}/\tilde b_{k+1}}{\tilde a_{k+1}/\tilde b_{k+1}-\tilde a_{k}/\tilde b_{k}}\text{ exists.}$$

\noindent{\bf(B2)$_{c}$} $\exists k_0>0,$ such that $\frac{\tilde a_k}{\tilde b_k}\ne\frac{\tilde a_{k+1}}{\tilde b_{k+1}},
\ \frac{\tilde d_k}{\tilde b_k}\ne\frac{\tilde d_{k+1}}{\tilde b_{k+1}},\ \forall k\ge k_0$ and
   $$\tau:=\lim_{k\rto}\frac{\tilde d_{k+1}/\tilde b_{k+1}-\tilde d_{k}/\tilde b_k}{\tilde a_{k+1}/\tilde b_{k+1}-\tilde a_{k}/\tilde b_k}\ne \frac{-(a+\theta)\pm \sqrt{(a+\theta)^2+4(bd-a\theta)}}{2b}  $$ exists as a finite or infinite number. In addition, if $\tau$ is finite,  assume further that $\lim_{k\rto}\frac{\tilde a_{k+2}/\tilde b_{k+2}-\tilde a_{k+1}/\tilde b_{k+1}}{\tilde a_{k+1}/\tilde b_{k+1}-\tilde a_{k}/\tilde b_{k}}
$ exists. Otherwise, if $\tau=\infty,$ assume  further that $\lim_{k\rto}\frac{\tilde d_{k+2}/\tilde b_{k+2}-\tilde d_{k+1}/\tilde b_{k+1}}{\tilde d_{k+1}/\tilde b_{k+1}-\tilde d_{k}/\tilde b_{k}}$ exists.

\begin{remark}\textbf{(i)} We remark that \eqref{abs} of condition (B1) is the requirement of \cite[Theorem 1]{hs20}, which we use to prove Theorem \ref{tmc} below. For all lemmas involved in this paper, we do not need \eqref{abs}. \textbf{(ii)} The conditions (B2)$_a,$ (B2)$_b$ and  (B2)$_c$ look a bit complicated. The lemma below gives an example for which (B1) and one of (B2)$_a,$ (B2)$_b$ and  (B2)$_c$ hold. Further example can also be found in Section \ref{eg}.
\end{remark}
\begin{lemma}\label{egc}
  Suppose that  $b,d>0$ and $a,\theta\ge0$ are numbers which are not all equal, satisfying $a+\theta>0,$ $\tau:=\frac{b(b+d-a-\theta)+2(a\theta-bd)}{b(2b-a-\theta)}\ne \frac{-(a+\theta)\pm\sqrt{(a+\theta)^2+4(bd-a\theta)}}{2b}$ and $(b-a)(b-\theta)\ge0.$ Let $r_n,n\ge 1$ be strictly positive numbers such that $\lim_{n\rto}r_n=0$ and $\lim_{n\rto}\frac{r_n-r_{n+1}}{r_n^{2}}=c$ for some number $0<c<\infty.$ Set $a_k=a+r_k,$ $b_k=b+r_k,$ $d_k=d+r_k$ and $\theta_k=\theta+r_k,$ $k\ge1.$ Then $a_k,b_k,d_k,\theta_k,k\ge1$ satisfies condition {(B1)} and one of the conditions {(B2)$_{ a},$ (B2)$_{ b}$} and {(B2)$_{ c}.$}
\end{lemma}
The proof of Lemma \ref{egc} is postponed to an Appendix.
Note that under condition (B1), we have \begin{equation}\label{lma}
  \lim_{k\rto}M_k=M:=\z(\begin{array}{cc}
                          a & b \\
                          d & \theta
                        \end{array}
  \y),\ \lim_{k\rto} A_k=A:=\z(\begin{array}{cc}
                          a+\theta & b \\
                          d-a\theta/b & 0
                        \end{array}
  \y)
\end{equation}
whose eigenvalues are  \begin{align*}
  \varrho(M)&=\varrho(A)=\frac{a+\theta+\sqrt{(a+\theta)^2+4(bd- a\theta)}}{2},\\
\varrho_1(M)&=\varrho_1(A)=\frac{a+\theta-\sqrt{(a+\theta)^2+4(bd-a\theta)}}{2}.
\end{align*}
Clearly, we have $|\varrho_1(A)|<\varrho(A).$ In the literature, the top eigenvalue $\varrho(A)$ is usually referred to as the spectral radius of $A.$
In what follows, we always denote by $\varrho(A)$  the spectral radius of a matrix $A.$

In the rest of the paper,  $f(n)\sim g(n)$ means $\lim_{n\rto}f(n)/g(n)=1,$ $f(n)=o(g(n))$ means $\lim_{n\rto}f(n)/g(n)=0$
and
 unless otherwise specified, $0<c<\infty$ is some constant, which may change from line to line. We always assume that empty product equals identity and empty sum equals $0.$ Now we are ready to state the main result.
 \begin{theorem}\label{pnec}
   Suppose that condition (B1) and one of the conditions (B2)$_{ a},$ (B2)$_{ b}$ and (B2)$_{ c}$ hold.
   Assume further  $|\varrho_1(M)|<1$  and  $\tilde d_k\ge\varepsilon$ for some $\varepsilon>0.$  Then
   \begin{align}\label{pns}
     &P(\nu>n|Z_0=\mb e_1)\sim \frac{c}{\sum_{k=1}^{n+1} \varrho(M_1)^{-1}\cdots \varrho(M_{k-1})^{-1}} \text{ as }n\rto.
   \end{align}
   If we assume $\theta\ne b+1$ in addition, then as $n\rto,$
   \begin{align}\label{pnse}
     &P(\nu=n|Z_0=\mb e_1)\sim\frac{c\varrho(M_1)^{-1}\cdots \varrho(M_n)^{-1}}{\z(\sum_{k=1}^{n+1} \varrho(M_1)^{-1}\cdots \varrho(M_{k-1})^{-1}\y)^2};
   \end{align}
   otherwise, if we suppose $\theta= b+1$ further, then as $n\rto,$
    \begin{align}\label{pnso}
     &P(\nu=n|Z_0=\mb e_1)=o\z(\frac{\varrho(M_1)^{-1}\cdots \varrho(M_n)^{-1}}{\big(\sum_{k=1}^{n+1} \varrho(M_1)^{-1}\cdots \varrho(M_{k-1})^{-1}\big)^2}\y).
   \end{align}
 \end{theorem}
 \begin{remark}
   \textbf{(a)} The righthand sides of \eqref{pns} and \eqref{pnse} are computable and easy to be estimated. In Section \ref{eg}, we give a example for which conditions of Theorem \ref{pnec} hold and $P(\nu=n|Z_0=\mb e_1)$ decays with various speeds, for example $\frac{c}{n(\log n)^2},$ $\frac{c}{n^\beta},\beta>1$ et al. as $n\rto,$ which are very different from those of homogeneous Galton-Watson processes. %
 \textbf{(b)} The condition $\forall k\ge1, \tilde d_k=d_k-a_k\theta_k/b_k\ge \varepsilon$ for some $\varepsilon>0$ (implying $bd>a\theta$) is a technical one, which ensures that  $A_k,k\ge1$ are  nonnegative matrices.
   If  $\forall k\ge1, \tilde d_k= d_k-a_k\theta_k/b_k\le -\varepsilon$ for some $\varepsilon>0$ (implying  $bd<a\theta$), $A_k,k\ge1$ are not nonnegative matrices any longer. But, in this case, by refining the proof of \cite[Theorem 1]{hs20}, at least when $\varrho(M)<1,$ our method does work.
    If $bd=a\theta,$ that is, the matrix $M$ is of rank one, our method does not work.
 \end{remark}

 \proof For the convenience of the reader, we give here the main body of the proof of Theorem \ref{pnec}. But the proofs of the four lemmas below,  especially those of Lemma \ref{gl} and Lemma \ref{smr}, are  the main difficulties of this paper and will be postponed to Section \ref{pbp} and Section \ref{pmp} below.

 To begin with, write $\tilde p_n\equiv P(\nu=n|Z_0=\mb e_1)$  and $ \tilde\eta_n\equiv P(\nu>n|Z_0=\mb e_1).$
 The lemma below gives explicitly the formulae of $\tilde p_n$ and $\tilde \eta_n.$
 \begin{lemma}\label{dna} For $n\ge1,$ we have
 \begin{align}\label{etan}
    \tilde\eta_n&=\frac{\mb e_1\prod_{k=1}^{n}A_k(1,\lambda_{n+1})^t}{\sum_{k=1}^{n+1}\mb e_1\prod_{i=k}^{n}A_i(1,\lambda_{n+1})^t},\\
\label{pn}
  \tilde p_n&=\frac{1}{\sum_{k=1}^{n+1}\mb e_1\prod_{i=k}^{n}A_i(1,\lambda_{n+1})^t}\frac{\mb e_1\prod_{k=1}^{n-1}A_{k}\mathbf e_1^t}{\sum_{k=1}^{n}\mb e_1\prod_{i=k}^{n-1}A_i(1,\lambda_{n})^t}G_{n-1}
  \end{align}
 where  $\lambda_n:=1-\frac{\theta_n}{b_n}$ and
\begin{align}\label{dg}
  G_{n-1}\equiv \frac{\mb e_1\prod_{k=1}^{n-1}A_k(1,\lambda_{n})^t\sum_{k=1}^{n+1}\mb e_1\prod_{i=k}^{n}A_i(1,\lambda_{n+1})^t}{\mb e_1\prod_{k=1}^{n-1}A_{k}\mathbf e_1^t}\\
  -\frac{\mb e_1\prod_{k=1}^{n}A_{k}(1,\lambda_{n+1})^t\sum_{k=1}^{n}\mb e_1\prod_{i=k}^{n-1}A_i(1,\lambda_{n})^t}{\mb e_1\prod_{k=1}^{n-1}A_{k}\mathbf e_1^t}.\no
  \end{align}
 \end{lemma}
 $G_n$ defined in \eqref{dg} looks very complicated. Its convergence in the next lemma plays a key role.
 \begin{lemma}\label{gl}
  Suppose that condition (B1) holds, $|\varrho_1(A)|<1$ and $\forall k, \tilde d_k\ge\varepsilon$ for some $\varepsilon>0.$  Then, $$
  \lim_{n\rto}G_n=G$$ exists. Moreover, we have $0<G<\infty$ if $\theta\ne b+1$ and $G=0$ if $\theta=b+1.$

\end{lemma}
 For $n\ge1,$ let $S_n:=\frac{\sum_{k=1}^{n+1}\varrho(A_1)^{-1}\cdots\varrho(A_{k-1})^{-1}}{\varrho(A_1)^{-1}\cdots\varrho(A_n)^{-1}}.
 $
 \begin{lemma} \label{smr} Under the conditions of Theorem \ref{pnec}, we have
 \begin{align}\label{epar}
   \mb e_1\prod_{i=1}^{n}A_n\mb e_1^t\sim c\varrho(A_1)\cdots\varrho(A_{n})\text{ as }n\rto
 \end{align} and for some number $0<\phi_1,\phi_2<\infty,$
   \begin{align}\label{asal}
    \lim_{n\rto}&\frac{\sum_{k=1}^{n+1}\mb e_1\prod_{i=k}^{n}A_i(1,\lambda_{n+1})^t}{S_n}=\phi_1 \text{ and } \lim_{n\rto}\frac{S_{n-1}}{S_n}=\phi_2.
  \end{align}
 \end{lemma}

 We only give here the proof of \eqref{pnse}, since \eqref{pns} and \eqref{pnso} can be proved similarly.
 Suppose now $\theta\ne b+1.$ Since by \eqref{lma}, we have $\lim_{n\rto} \varrho(A_n)=\varrho(A)>0.$   Thus it follows immediately from the above three lemmas that
 \begin{align}\label{pna}
     &\tilde p_n\sim\frac{c\varrho(A_1)^{-1}\cdots \varrho(A_n)^{-1}}{\z(\sum_{k=1}^{n+1} \varrho(A_1)^{-1}\cdots \varrho(A_{k-1})^{-1}\y)^2}.
   \end{align}
 The lemma below allows us to go from the matrices $A_k$ back to $M_k,$ $k\ge1.$
 \begin{lemma}\label{bma}
 Suppose the condition (B1) is satisfied and $\tilde d_k\ge \varepsilon,\forall k\ge1$ for some $\varepsilon>0.$ Then as $n\rto,$
   $$\frac{\varrho(A_1)^{-1}\cdots \varrho(A_n)^{-1}}{\z(\sum_{k=1}^{n+1} \varrho(A_1)^{-1}\cdots \varrho(A_{k-1})^{-1}\y)^2}\sim \frac{c\varrho(M_1)^{-1}\cdots \varrho(M_n)^{-1}}{\z(\sum_{k=1}^{n+1} \varrho(M_1)^{-1}\cdots \varrho(M_{k-1})^{-1}\y)^2}.$$
  \end{lemma}
 As a result, taking \eqref{pna} and Lemma \ref{bma} together, we get \eqref{pnse}.  \qed

%

 Lemma \ref{smr} is a consequence of Proposition \ref{tmc} below, which may have its own interest in the theory of product of positive matrices.
 Let \begin{align*}
   B_k=\left(
         \begin{array}{cc}
           a_k & b_k \\
           d_k & 0 \\
         \end{array}
       \right),k\ge1.
 \end{align*}
 \begin{proposition}\label{tmc}
      Suppose that $a_k,b_k,d_k,\theta_k,$ $k\ge1$ satisfies  condition (B1) and one of the conditions (B2)$_a,$ (B2)$_b$ and  (B2)$_c,$ with $\theta_k\equiv 0,k\ge1$ and the number $a>0.$ Then  $\forall i,j\in \{1,2\},$ we have
  \begin{align}\label{rac}
    \mathbf e_iB_1\cdots B_n\mathbf e_j^t\sim c\varrho(B_1)\cdots\varrho(B_n), \text{ as }n\rto.
  \end{align}Furthermore, writing $$\tilde S_n:=\frac{\sum_{k=1}^{n+1}\varrho(B_1)^{-1}\cdots\varrho(B_{k-1})^{-1}}{\varrho(B_1)^{-1}\cdots\varrho(B_n)^{-1}},\ Y_n:=\frac{\sum_{k=1}^{n+1}\mb e_1B_k\cdots B_n\mb e_1^t}{\tilde S_n},n\ge0,$$
  then for some number $\psi>0,$ we have
 \begin{align}\label{sax}
   \lim_{n\rto}Y_n=\psi.
  \end{align}
 \end{proposition}
\begin{remark} {\bf(i)} The limit number $\psi$ can be deduced from \eqref{ttt} and \eqref{psi} below.
  {\bf(ii)} Even though we have known \eqref{rac}, it is hard to evaluate $\sum_{k=1}^{n+1}\mb e_1B_k\cdots B_n\mb e_1^t$ since every summand there depends on $n.$ But by \eqref{sax}, $\sum_{k=1}^{n+1}\mb e_1B_k\cdots B_n\mb e_1^t $ is asymptotically the same as $\tilde S_n,$ which is computable.
\end{remark}

 \section{Examples}\label{eg}
For $n\ge1,$ let $q_{n,1}\ge0,q_{n,2}>0, p_n>0$ be numbers such that $q_{n,1}+q_{n,2}+p_n=1.$ Suppose that $Z_n, n\ge0$ is a 2-type branching process with $Z_0=\mathbf e_1$ and offspring distributions
\begin{align}
\label{pha}&P(Z_{n}=(i,j)|Z_{n-1}=\mb e_{1})=\frac{(i+j)!}{i!j!}q_{n,1}^{i}q_{n,2}^{j}p_{n},\\
\label{phb}&P(Z_{n}=(1+i,j)|Z_{n-1}=\mb e_{2})=\frac{(i+j)!}{i!j!}q_{n,1}^{i}q_{n,2}^{j}p_{n},
\end{align}where $i\ge0,j\ge 0,n\ge1.$
Some computation yields that
the mean matrix
\begin{align*}
M_n:=\z(\begin{array}{cc}
          a_n & b_n \\
          1+a_n & b_n
        \end{array}
\y)\text{ with }a_n:=\frac{q_{n,1}}{p_{n}},b_n:=\frac{q_{n,2}}{p_{n}},n\ge1.
\end{align*}
Also, it follows from \eqref{pha} and \eqref{phb} that
\begin{align}
F_{n}(\mathbf{s})&:=E(s_{1}^{Z_{n,1}}s_{2}^{Z_{n,2}}|Z_0=\mb e_1)
=1-\frac{\mathbf{e}_{1}\prod_{j=1}^{n}M_{j}(\mb 1-\mathbf{s})}{1+\sum_{k=1}^{n}\mathbf{e}_{1}\prod_{i=k}^{n}M_{i}(\mb1-\mathbf{s})}\no
\end{align}
which leads to
\begin{align}
  \label{pnh}\eta_n:&=P(\nu=n|Z_0=\mb e_1)=F_n(\mb 0)-F_{n-1}(\mb 0)\\
  &=\frac{\mathbf{e}_{1}\prod_{j=1}^{n-1}
  M_{j}\mb1}{\sum_{k=1}^{n}\mathbf{e}_{1}\prod_{i=k}^{n-1}M_{i}\mb1}-\frac{\mathbf{e}_{1}\prod_{j=1}^{n}M_{j}\mb 1}{\sum_{k=1}^{n+1}\mathbf{e}_{1}\prod_{i=k}^{n}M_{i}\mb1}.\no
\end{align}
At first, we see what happens to the homogeneous case, that is \begin{align}
  q_{n,1}\equiv q_1,q_{n,2}\equiv q_2,p_n\equiv p \text{ so that }M_n\equiv M:=\z(\begin{array}{cc}
          q_1/p & q_2/p \\
          1+q_1/p & q_2/p
        \end{array}
\y),n\ge1,\no
\end{align}
where we assume $p>0,q_1\ge0, q_2>0$ and $p+q_1+q_2=1$ which ensure that $Z_n,n\ge0$ is truly a two-type branching process.
Clearly,  the spectral radius of $M$ is $$\varrho(M)=\frac{q_1+q_2+\sqrt{(q_1+q_2)^2+4pq_2}}{2p}.$$
By some easy computation, from \eqref{pnh}, we get that as $n\rto,$
\begin{align*}
  \eta_n&=\frac{\mb e_1M^{n-1}\mb 1}{\sum_{k=0}^{n-1}\mb e_1M^k\mb 1}-\frac{\mb e_1M^n\mb 1}{\sum_{k=0}^n\mb e_1M^k\mb 1}
  \sim \left\{\begin{array}{ll}
                 c\varrho(M)^n, &\text{if }\varrho(M)<1, \\
                 c\varrho(M)^{-n}, &\text{if }\varrho(M)>1, \\
                 n^{-2}, &\text{if }\varrho(M)=1.
               \end{array}
  \right.
\end{align*}
We thus come to the conclusion that for homogeneous Galton-Watson processes,
$P(\nu=n|Z_0=\mb e_1)$ decays either {\it exponentially} (supercritical and subcritical cases) or {\it polynomially} (critical case).

Next, we turn to consider BPVE by adding some perturbations on a critical homogeneous Galton-Watson process, which may exhibit very different asymptotics.
For $K= 1,2,...$ and $B\in \mathbb R,$ set
\begin{align*}
    &\Lambda(1,i,B)=\frac{B}{i},\\
   &\Lambda(2,i,B)=\frac{1}{i}+\frac{B}{i\log i}, \cdots,\\
     &\Lambda(K,i,B)=\frac{1}{i}+\frac{1}{i\log i}+...+\frac{1}{i\log i\cdots\log_{K-2}i}+\frac{B}{i\log i\cdots\log_{K-1}i},
  \end{align*}
  where $\log_0 i=i$ and for $k\ge1,$ $\log_{k}i=\log \log_{k-1} i.$
Now fix $K$ and $B,$  set $i_0:=\min\z\{i: \log_{K-1}i>0, {|\Lambda(K,i,B)|}< 1 \y\}$
 and let
\begin{align*}r_i:=\left\{\begin{array}{ll}
   \frac{\Lambda(K,i,B)}{3}, & i\ge  i_0, \\
 r_{i_0}, & i< i_0,
\end{array}\right.\end{align*}
which serve as perturbations. To avoid tedious computation, we assume $q_{n,1}\equiv 0,\forall n\ge1$ and write $q_{n,2}$ simply as $q_n$ for $n\ge1$ so that \eqref{pha} and \eqref{phb} reduce to
\begin{align}
  \label{phaa}&P(Z_{n}=(0,j)|Z_{n-1}=\mb e_{1})=q_{n}^{j}p_{n},\\
\label{phbb}&P(Z_{n}=(1,j)|Z_{n-1}=\mb e_{2})=q_{n}^{j}p_{n},j\ge0,n\ge1.
\end{align}

\begin{theorem}\label{mxt} Suppose that $Z_n,n\ge0$ is a two-type branching process whose offspring distribution satisfies \eqref{phaa} and \eqref{phbb}.  Fix $K=1,2,3,...$ and $B\in \mathbb R.$
 {\rm (i)} If
    $p_i=\frac{2}{3}+r_i,q_i=\frac{1}{3}-r_i, i\ge 1,$
then, as $n\rto,$
\begin{align*}
  P(\nu=n|Z_0=\mb e_1)\sim \left\{ \begin{array}{ll}
                                \frac{c }{n\log n\cdots \log_{K-2}n\log_{K-1} n(\log_K n)^2}, &  \text{if } B=1,  \\
\frac{c}{n\log n\cdots \log_{K-2}n(\log_{K-1} n)^B}, &  \text{if } B>1,  \\
\frac{c}{n\log n\cdots \log_{K-2}n(\log_{K-1} n)^{2-B}}, &  \text{if } B<1.
                              \end{array}
  \right.
\end{align*}
{\rm(ii)} If
$  p_i=\frac{2}{3}-r_i,q_i=\frac{1}{3}+r_i,i\ge 1,$
then, as $n\rto,$
\begin{align*}
  P(\nu=n|Z_0=\mb e_1)\sim
  \left\{\begin{array}{ll}
\frac{c}{n^{B+2}}, &\text{if } K=1,B>-1,\\
\frac{c}{n(\log n)^2}, & \text{if }K=1, B=-1,\\
  cn^B, & \text{if }K=1,B<-1,\\
 \frac{c}{n^3\log n...\log_{K-2}n (\log_{K-1} n)^B}, &\text{if } K>1.\\
\end{array}\right.
\end{align*}
\end{theorem}
\begin{remark}
  As seen from above, for two-type homogeneous  processes, $P(\nu=n|Z_0=\mb e_1)$ decays either exponentially or polynomially with speed $n^{-2}$. But for two-type BPVE, $P(\nu=n|Z_0=\mb e_1)$ may decay with many different speeds such as $\frac{c}{(n\log n)^2},$ $cn^{-B},  B>1$  et al. as $n\rto.$
\end{remark}
\proof Write $b_k=q_k/p_k,k\ge1.$ Then $M_k=\left(
                  \begin{array}{cc}
                 0 & b_k \\
                    1 & b_k \\
 \end{array}
\right)
 \rightarrow M=\left(
        \begin{array}{cc}
             0 & 1/2 \\
             1 & 1/2 \\
             \end{array}
     \right)$ as $k\rto.$ Comparing $M_k$ with the one in \eqref{mg}, we find that $a_k\equiv0,$ $\theta_k=b_k,$ $d_k\equiv 1,k\ge1$ and  $\varrho(M)=1,$ $\varrho_1(M)=-1/2.$ By definition, we know that $\tilde a_k=b_k,$ $\tilde b_k=b_k,$ $\tilde d_k=d_k\equiv 1.$  Thus, $A_k=\left(
        \begin{array}{cc}
             b_k & b_k \\
             1 & 0 \\
             \end{array}
     \right)$ and it is easily seen that $\varrho(A_k)\equiv\varrho(M_k),k\ge1.$
     Moreover, we have the following lemma.
\begin{lemma}\label{dfr}{\bf (i)} $\lim\limits_{n\rightarrow\infty}\frac{r_n-r_{n+1}}{n^2}=1/3$ and thus $\sum_{k=1}^\infty|b_{k+1}-b_k|<\infty$ whenever $p_i=2/3\pm r_i, i\ge 1.$
{\bf (ii)}  We have $\frac{1}{b_k}\neq\frac{1}{b_{k+1}},$ $\forall k\ge i_0$ and
$\frac{b_{k+1}-b_k}{b_{k+2}-b_{k+1}}\rightarrow1$ as $n\rto,$ no mater $p_i=2/3+r_i$ or $p_i=2/3-r_i,i\ge 1.$
\end{lemma}
For the proof of the lemma, we refer the reader to \cite[Lemma 7]{hs20}. It follows immediately from Lemma \ref{dfr} that condition (B1) and condition  (B2)$_{a}$ are fulfilled.
Therefore, all conditions of Theorem \ref{pnec} are satisfied. Consequently, applying Theorem \ref{pnec}, we conclude that
\begin{align}\label{pnsa}
     P(\nu=n|Z_0=\mb e_1)\sim \frac{c\varrho(M_1)^{-1}\cdots \varrho(M_n)^{-1}}{\z(\sum_{k=1}^{n+1} \varrho(M_1)^{-1}\cdots \varrho(M_{k-1})^{-1}\y)^2}, \text{ as }n\rto.
   \end{align}

  Note that $\varrho(M_k)=\z(b_k+\sqrt{b_k^2+4b_k}\y)/2.$ If $p_i=\frac{2}{3}\pm r_i,i\ge 1,$ then by Taylor enpension of  $\varrho(M_k)$ at $0,$ we get
  \begin{equation*}\varrho(M_k)=1\mp 3r_k+O(r_k^2)\text{ as }k\rto.\end{equation*}
Applying \cite[Proposition 2]{hs20}, we get
\begin{align}\varrho(M_1)\cdots\varrho(M_n)\sim c\z(n\log n\cdots \log_{K-2}n(\log_{K-1} n)^B\y)^{\mp 1}.\label{rpd}\end{align}
With \eqref{pnsa} and \eqref{rpd} in hands, the proof of Theorem \ref{mxt} goes almost verbatim as that of \cite[Theorem 1]{wh}. We do not repeat it here. \qed
 \section{Proofs of the auxiliary Lemmas }\label{pbp}
  The main task of this section is to finish the proofs of Lemmas \ref{dna}, \ref{gl}, \ref{smr} and \ref{bma}, which are required when proving Theorem \ref{pnec}. The proof of Lemma \ref{smr} is based on Proposition \ref{tmc} whose proof is very long and will be postponed to Section \ref{pmp}.
\subsection{Proof of Lemma \ref{dna}}\label{ta}

For $n\ge 0$ and $\mathbf s=(s_1,s_2)^t\in [0,1]^2,$  let $$F_n^{(i)}(\mathbf s)\equiv E(\mathbf s^{Z_n}|Z_0=\mathbf e_i):=E( s_1^{Z_{n,1}}s_2^{Z_{n,2}}|Z_0=\mathbf e_i),i=1,2$$
  and set $$ \mathbf F_n(\mathbf s)=(F_n^{(1)}(\mathbf s),F_n^{(2)}(\mathbf s))^t.$$
It follows by induction (see Dyakonova \cite[Lemma 1]{dy99}) that for $n\ge0,$
\begin{align*}
  \mathbf F_{n}(\mathbf{s})=\mathbf{f}_{1}(\mathbf{f}_{2}(\cdots\mathbf{f}_{n}(\mathbf{s})\cdots))
  =\mathbf1-\frac{\prod_{k=1}^{n}M_{k}(\mathbf 1-\mathbf{s})}{1+\sum_{k=1}^{n}\gamma_k\prod_{i=k+1}^{n}M_{i}(\mathbf1-\mathbf{s})}
\end{align*}
which leads to
\begin{align*}
  \mathbf F_{n}(\mathbf{0})=\mathbf 1- \frac{\prod_{k=1}^{n}M_{k}\mathbf 1}{1+\sum_{k=1}^{n}\gamma_k\prod_{i=k+1}^{n}M_{i}\mathbf1}.
\end{align*}
 Note that
for $n\ge1$ and $i=1,2,$
\begin{align}\label{pngn}
  P(\nu>n|Z_0=\mb e_i)&=1-F_n^{(i)}(\mb 0)=\frac{\mb e_i\prod_{k=1}^{n}M_{k}\mathbf 1}{\sum_{k=1}^{n+1}\gamma_k\prod_{i=k+1}^{n}M_{i}\mathbf1},
\end{align}
and consequently \begin{align}\label{pnen}
  P(&\nu=n|Z_0=\mb e_i)=\frac{\mb e_i\prod_{k=1}^{n-1}M_{k}\mathbf 1}{\sum_{k=1}^{n}\gamma_k\prod_{i=k+1}^{n-1}M_{i}\mathbf1}-\frac{\mb e_i\prod_{k=1}^{n}M_{k}\mathbf 1}{\sum_{k=1}^{n+1}\gamma_k\prod_{i=k+1}^{n}M_{i}\mathbf1}.
\end{align}
%
With matrices  $A_i,i\ge1$ the ones defined in \eqref{dta},
   using \eqref{mg} and \eqref{am},  we have from \eqref{pngn} and \eqref{pnen} that \eqref{etan} and \eqref{pn} hold. Thus Lemma \ref{dna} is proved. \qed
%

\subsection{Proof of Lemma \ref{gl}}
Due to the complicated formula of $G_n,$ the proof will be a long journey. So we divide the proof into several steps. Keep in mind that $\lambda_n:=1-\theta_n/b_n\rightarrow 1-\theta/b$ as $n\rto.$

{\it Step 1: We show that $\lim_{n\rto}G_n=G$ for some number $0\le G<\infty.$} For this purpose,
 set  \begin{align*}
  f_n\equiv \frac{\mb e_1\prod_{k=1}^{n}A_{k}\mathbf e_2^t}{\mb e_1\prod_{k=1}^{n}A_{k}\mathbf e_1^t} \text{ and }H_n\equiv\sum_{k=1}^{n}\mb e_1\prod_{i=k}^{n}A_{i}(f_n\mathbf e_1^t-\mb e_2^t),n\ge1.
\end{align*}
Then by some subtle computation, we get
\begin{align}\label{ghf}
  G_{n-1}=1&+(\tilde b_n\lambda_n\lambda_{n+1} +\tilde a_n\lambda_{n}-\tilde d_n)H_{n-1}\\
  &+(\tilde b_n\lambda_n\lambda_{n+1} +\tilde a_n\lambda_{n}-\tilde d_n+\lambda_{n})f_{n-1}.\no
\end{align}
The lemma below, whose proof is postponed to the end of this subsection,  confirms the convergence of $f_n$ and $H_n$ as $n\rto.$
\begin{lemma}\label{hs}
  We have
  \begin{align}
    H_n=\sum_{k=1}^{n-1}(-1)^{n-k}f_k\tilde d_{k+1}f_{k+1}\cdots \tilde d_nf_n,n\ge1.\no
  \end{align}
 Furthermore, if condition (B1) holds, $|\varrho_1(A)|<1$ and $\forall k \ge1, \tilde d_k\ge \varepsilon$ for some $\varepsilon>0,$   then $$\lim_{n\rto}f_n=-\frac{b\varrho_1(A)}{bd- a\theta}\text{ and } \lim_{n\rto}H_n=-\frac{b}{bd- a\theta}\frac{\varrho_1(A)^2}{1-\varrho_1(A)}.$$
\end{lemma}

Applying Lemma \ref{hs}, from \eqref{ghf} we get
\begin{align*}
    \lim_{n\rto}G_n=G:=\frac{(b-\theta)\varrho_1(A)^2-(b-\theta)(a+b+1)\varrho_1(A)+bd-a\theta }{(bd-a\theta)(1-\varrho_1(A))}.
  \end{align*}

{\it Step 2: Show that $0<G<\infty$ if $\theta\ne b+1$ and $G=0$ if $\theta=b+1.$}
  To begin with, we prove that  $G$ is nonnegative. Indeed, in view of \eqref{pn}, since  $\sum_{k=1}^{n}\mb e_1\prod_{i=k}^{n-1}A_i(1,\lambda_{n})^t=\sum_{k=1}^{n}\mb e_1\prod_{i=k}^{n-1}M_{i}\mathbf1>0$ and $\mb e_1\prod_{k=1}^{n-1}A_{k}\mathbf e_1^t=\mb e_1\prod_{k=1}^{n-1}M_{k}\mathbf (1,\theta_{k+1}/b_{k+1})^t>0,$ we must have $G_n\ge 0,\forall n\ge0$ and hence $G=\lim_{n\rto}G_n\ge 0.$

  Note that by assumption, we have $\forall k\ge1,\tilde d_k\ge\varepsilon$ for some $\varepsilon>0,$ so we must have $bd> a\theta$ and hence $\varrho_1(M)<0.$ If $b\ge \theta,$ it is easily seen that $G>0.$ Thus we suppose next $b<\theta.$
  %
%
Let $$g(x)=(b-\theta)x^2-(b-\theta)(a+b+1)x+bd- a\theta.$$ Then $g(x)=0$ has two roots
$$\frac{1}{2}\z(a+b+1\pm\sqrt{(a+b+1)^2+4\frac{bd- a\theta}{\theta-b}}\y).$$
Since $\varrho_1(M)<0,$ then
$G=0$ if and only if $$\varrho_1(A)=\frac{1}{2}\z(a+b+1-\sqrt{(a+b+1)^2+4\frac{bd- a\theta}{\theta-b}}\y),$$
or equivalently,
\begin{align}\label{rx}
  a+\theta+&\sqrt{(a+ b+1)^2+4\frac{bd-a\theta}{\theta-b}}\\
  &=(a+b+1)+\sqrt{(a+\theta)^2+4(bd-a\theta)}.\no
\end{align}
By some subtle computation, we see that \eqref{rx} happens if and only if
\begin{align*}
  &(bd-a\theta)\z(b+1-\theta\y)\\
  &\quad\times\bigg\{2(a^2+b^2+2ab+b-\theta)+(a+b+1)\sqrt{(a+\theta)^2+4(b d-a\theta)} \\
  &\quad\quad\quad\quad+ (a+\theta)\sqrt{(a+ b+1)^2+4\frac{bd-a\theta}{\theta-b}} \bigg\}\\
  &=:(bd-a\theta)\z(b+1-\theta\y)\Theta(a,b,d,\theta)=0.
\end{align*}
Clearly, we have $\Theta(a,b,d,\theta)>0.$ Since $bd>a\theta,$  we come to the conclusion that $G=0$ if and only if $\theta = b+1.$

Consequently, if Lemma \ref{hs} is true, taking the above two steps together, we conclude that Lemma \ref{gl} holds.
\qed

To end this subsection, we prove Lemma \ref{hs}.

\noindent{\it Proof of Lemma \ref{hs}.} To begin with, we prove the first part. Note that
\begin{align}
  \label{fr}f_n=\frac{\mb e_1A_1\cdots A_{n}\mathbf e_2^t}{\mb e_1A_{1}\cdots A_n\mathbf e_1^t}=\frac{\tilde b_n\mb e_1A_1\cdots A_{n-1}\mb e_1^t}{\mb e_1A_{1}\cdots A_{n-1}(\tilde a_n\mathbf e_1^t+\tilde d_n\mb e_2^t)}=\frac{\tilde b_n}{\tilde a_n+\tilde d_nf_{n-1}},
\end{align}
which leads to
\begin{align}
  \tilde d_nf_nf_{n-1}=\tilde b_n-\tilde a_nf_{n},n\ge2.\no
\end{align}
Consequently, for $n\ge2,$
\begin{align}
  H_n&=\sum_{k=1}^{n}\mb e_1\prod_{i=k}^{n}A_{i}(f_n\mathbf e_1^t-\mb e_2^t)\label{hr}\\
  &=\tilde a_nf_n-\tilde b_n+\sum_{k=1}^{n-1}\mb e_1\prod_{i=k}^{n-1}A_{i}((\tilde a_nf_n-\tilde b_n)\mathbf e_1^t+\tilde d_nf_n\mb e_2^t)\no\\
  &=- \tilde d_nf_nf_{n-1} -\tilde d_nf_n\sum_{k=1}^{n-1}\mb e_1\prod_{i=k}^{n-1}A_{i}(f_{n-1}\mathbf e_1^t-\mb e_2^t)\no\\
  &=- \tilde d_nf_nf_{n-1} -\tilde d_nf_n H_{n-1}.\no
\end{align}
Since $H_1=0,$ iterating \eqref{hr}, we get
$$H_n=\sum_{k=1}^{n-1}(-1)^{n-k}f_k\tilde d_{k+1}f_{k+1}\cdots \tilde d_nf_n,n\ge1.$$
Next, we turn to prove the second part.  Suppose that condition (B1) holds, $|\varrho_1(M)|<1$ and $\forall k, \tilde d_k\ge \varepsilon$ for some $\varepsilon>0.$ Then it is easily seen that $bd> a\theta.$ Iterating \eqref{fr}, we get
\begin{align*}
  f_n=\frac{\tilde b_n\tilde d_n^{-1}}{\tilde a_n\tilde d_n^{-1}}\begin{array}{c}
                                \\
                               +
                             \end{array}\frac{\tilde b_{n-1}\tilde d_{n-1}^{-1}}{ \tilde a_{n-1}\tilde d_{n-1}^{-1}}\begin{array}{c}
                                \\
                               +\cdots+
                             \end{array}\frac{\tilde b_1\tilde d_1^{-1}}{\tilde a_1\tilde d_1^{-1}},n\ge1.
\end{align*}
 Thus by the theory of  convergence of limit periodic continued fractions(see Lorentzen and Waadeland \cite[Theorem 4.13, page 188]{lw}), the limit $f:=\lim_{n\rto}f_n$ exists. But since $ f_n>0, \forall n\ge1,
$ we must have $f\ge 0.$
Letting $n\rto$ in \eqref{fr}, we get
$f=\frac{b}{a+\theta +(d-a\theta/b)f}$ whose positive solution is
$$f= \frac{b\z(\sqrt{(a+\theta)^2+4(b d- a\theta)}-a-\theta\y)}{2(bd- a\theta)}=-\frac{b\varrho_1(A)}{bd-a\theta}.$$
Consequently,  we have $$\lim_{n\rto}\tilde d_nf_n=-\varrho_1(A)\in(0,1).$$
For any $\epsilon\in (0,1+\varrho_1(A)),$ there exists a number $k_1>0$ such that $- \varrho_1(A)-\epsilon< \tilde d_kf_k<-\varrho_1(A)+\epsilon, f-\epsilon<f_k<f+\epsilon,\forall k\ge k_1.$
Write
\begin{align}
 \label{hn}H_n&=\sum_{k=1}^{n-1}(-1)^{n-k}f_k\tilde d_{k+1}f_{k+1}\cdots \tilde d_nf_n\\
  & =\sum_{k=1}^{k_1-1}(-1)^{n-k}f_k\tilde d_{k+1}f_{k+1}\cdots \tilde d_nf_n+\sum_{k=k_1}^{n-1}(-1)^{n-k}f_k\tilde d_{k+1}f_{k+1}\cdots \tilde d_nf_n\no\\
  &=: H_n(1)+H_n(2).\no
\end{align}
Since $\lim_{n\rto}f_n=f,$ $f_n\le C, \forall n\ge1$ for some number $C>0.$ Thus we have
\begin{align}
  \label{lhna} |H_n(1)|\le C\sum_{k=1}^{k_1-1}(-\varrho_1(A))^{n-k_1+1}=C(k_1-1)(-\varrho_1(A))^{n-k_1+1}\rightarrow0
\end{align}
as $n\rto.$ On the other hand,
\begin{align}
  &\varlimsup_{n\rto}H_n(2)\le (f+\epsilon) \frac{(-\varrho_1(A)+\epsilon)^2}{1-(-\varrho_1(A)+\epsilon)^2}-(f-\epsilon)\frac{(-\varrho_1(A)-\epsilon)}{1-(-\varrho_1(A)-\epsilon)^2},\no\\
   &\varliminf_{n\rto}H_n(2)\ge (f-\epsilon) \frac{(-\varrho_1(A)-\epsilon)^2}{1-(-\varrho_1(A)-\epsilon)^2}-(f+\epsilon)\frac{(-\varrho_1(A)+\epsilon)}{1-(-\varrho_1(A)+\epsilon)^2}\no
\end{align} from which we get by letting $\epsilon\rightarrow 0$ that
\begin{align}\label{lhnb}
   \lim_{n\rto}H_n(2)= f\z( \frac{\varrho_1(A)^2}{1-\varrho_1(A)^2}+\frac{\varrho_1(A)}{1-\varrho_1(A)^2}\y)=\frac{f\varrho_1(A)}{1-\varrho_1(A)}.
 \end{align}
Taking limits on both sides of \eqref{hn}, we get from \eqref{lhna} and \eqref{lhnb} that
$$\lim_{n\rto}H_n=\frac{f\varrho_1(A)}{1-\varrho_1(A)}=-\frac{b}{bd-a\theta }\frac{\varrho_1(A)^2}{1-\varrho_1(A)}.$$
The lemma is proved.\qed

\subsection{Proof of Lemma \ref{smr}}
To begin with,  we prove the following lemma which  will be used times and again.
  \begin{lemma}\label{sig} Suppose that $\sigma_n,n\ge1$ is a sequence of positive numbers and $\lim_{n\rto}\sigma_n=\sigma>0.$ Then we have
    \begin{align}\label{sss}
    \lim_{n\rto}\frac{\sigma_1\cdots\sigma_{n+1}}{\sum_{k=1}^{n+1}\sigma_1\cdots\sigma_{k-1}}=\z\{\begin{array}{cc}
                                          0, & \text{if }\sigma\le1, \\
                                          \sigma-1, & \text{if }\sigma>1.
                                        \end{array}
    \y.
  \end{align}
  \end{lemma}
\proof
  Suppose first  $\sigma\le 1.$ Clearly, if $\sigma<1,$ then  $\lim_{n\rto}\sigma_1\cdots \sigma_n=0$ and \eqref{sss} holds trivially. Now we assume $\sigma=1.$ Fix $N>0.$ For $n>N,$ we have
   \begin{align*}
     &\frac{\sum_{k=1}^{n+1}\sigma_1\cdots\sigma_{k-1}}{\sigma_1\cdots\sigma_{n+1}}=\sum_{k=1}^{n+1}\sigma_k^{-1}\cdots\sigma_{n+1}^{-1}\\
     &\quad\quad\ge \sigma_{n+1}^{-1}\z( 1+\frac{1}{\sigma_n}+\frac{1}{\sigma_n\sigma_{n-1}}+\dots+\frac{1}{\sigma_n\cdots\sigma_{n-N+1}}\y).
   \end{align*}
   Letting $n\rto,$  we have $\varliminf_{n\rto}\frac{\sum_{k=1}^{n+1}\sigma_1\cdots\sigma_{k-1}}{\sigma_1\cdots\sigma_{n+1}}\ge (N+1)\sigma^{-1}.$ Since $N$ is arbitrary, we have $\lim_{n\rto}\frac{\sum_{k=1}^n\sigma_1\cdots\sigma_k}{\sigma_1\cdots\sigma_{n+1}}=\infty$ which implies \eqref{sss}.

      Suppose next $\sigma> 1.$ Since $\lim_{n\rto}\sigma_n=\sigma>0,$ for $0<\epsilon<(1-\sigma^{-1})\wedge \sigma^{-1},$ $\exists k_2>0$ such that $\sigma^{-1}-\epsilon\le \sigma_k^{-1}\le \sigma^{-1}+\epsilon,\forall k>k_2.$ With this number $k_2,$ we get
        \begin{align}\label{ds}
        \frac{\sum_{k=1}^{n+1}\sigma_1\cdots\sigma_{k-1}}{\sigma_1\cdots\sigma_{n+1}}=\sum_{k=1}^{k_2}\sigma_k^{-1}\cdots \sigma_{n+1}^{-1}+\sum_{k=k_2+1}^{n+1}\sigma_k^{-1}\cdots\sigma_{n+1}^{-1}.
      \end{align}
      Clearly we have $0\le \sum_{k=1}^{k_2}\sigma_k^{-1}\cdots \sigma_{n+1}^{-1} \le \sum_{k=1}^{k_2}\sigma_k^{-1}\cdots \sigma_{k_2}^{-1} (\sigma^{-1}+\epsilon)^{n-k_2+1}\rightarrow 0$ as $n\rto$ and
      \begin{align*}
        \frac{\sigma^{-1}-\epsilon}{1-(\sigma^{-1}-\epsilon)}&\le\varliminf_{n\rto}\sum_{k=k_2+1}^{n+1}\sigma_k^{-1}\cdots\sigma_{n+1}^{-1}\\
      &\le\varlimsup_{n\rto}\sum_{k=k_2+1}^{n+1}\sigma_k^{-1}\cdots\sigma_{n+1}^{-1}\le \frac{\sigma^{-1}+\epsilon}{1-(\sigma^{-1}+\epsilon)}.
      \end{align*} Letting $\epsilon\rightarrow0,$ we get
      $\lim_{n\rto}\sum_{k=k_2+1}^{n+1}\sigma_k^{-1}\cdots\sigma_{n+1}^{-1}=(\sigma-1)^{-1},$ which together with \eqref{ds} implies \eqref{sss}. \qed

 Now, suppose that condition (B1) holds. We show that
  \begin{align}\label{srl}
    \lim_{n\rto}\frac{S_{n-1}}{S_n}=\z\{\begin{array}{cc}
                                          \varrho(A)^{-1}, & \text{if }\varrho(A)\ge1, \\
                                          1, & \text{if }\varrho(A)<1.
                                        \end{array}
    \y.
  \end{align}
  Indeed, by some easy computation, we obtain
  \begin{align}\label{sd}
    \frac{S_{n-1}}{S_n}&=\varrho(A_n)^{-1}\frac{1}{1+\frac{\varrho(A_1)^{-1}\cdots\varrho(A_{n})^{-1}}{\sum_{k=1}^{n}\varrho(A_1)^{-1}\cdots\varrho(A_{k-1})^{-1}}}.
  \end{align}
  Note that $\varrho(A_n)>0,\forall n\ge1$ and  under condition (B1), $\lim_{n\rto}\varrho(A_n)=\varrho(A).$ Thus we can apply Lemma \ref{sig} to the sequence $\varrho(A_n)^{-1},n\ge1$ to get
  \begin{align*}
    \lim_{n\rto}\frac{\varrho(A_1)^{-1}\cdots\varrho(A_{n})^{-1}}{\sum_{k=1}^{n}\varrho(A_1)^{-1}\cdots\varrho(A_{k-1})^{-1}}=\z\{\begin{array}{cc}
                                          0, & \text{if }\varrho(A)\ge1, \\
                                          \varrho(A)^{-1}-1, & \text{if }\varrho(A)<1.
                                        \end{array}
    \y.
  \end{align*}
  Then taking limits on both sides of \eqref{sd}, we get \eqref{srl}.
  The second limit in \eqref{asal} is proved.

Note that under the assumptions of Theorem \ref{pnec}, we have  $\tilde a_k\rightarrow a+\theta>0,$ $\tilde b_k\rightarrow b>0,$ $\tilde d_k\rightarrow d-a\theta/b>0$ as $k\rto.$
   Moreover  it is easy to check that $\sum_{k=2}^\infty |\tilde a_{k}-\tilde a_{k-1}|+|\tilde b_{k}-\tilde b_{k-1}|+|\tilde d_{k}-\tilde d_{k-1}|<\infty.$
   By assumption,
   $\tilde a_k,\tilde b_k,\tilde d_k,k\ge1$ satisfies  one of the conditions (B2)$_{a},$ (B2)$_{b}$ and (B2)$_{c}.$ Therefore all conditions of Proposition \ref{tmc} are also fulfilled for $\tilde a_k,\tilde b_k, \tilde d_k,k\ge1.$
  Let \begin{align*}
  \tilde Y_n:=\frac{\sum_{k=1}^{n+1}\mb e_1A_k\cdots A_n\mb e_1^t}{S_n},n\ge0.\end{align*}
  Then applying Proposition \ref{tmc} to the matrices $A_k,k\ge1,$ we get from \eqref{rac} and  \eqref{sax} that, with some number $0<\tilde \psi<\infty,$ we have
 \begin{align}
    \label{tya}&\lim_{n\rto}\tilde Y_n=\tilde\psi,\text{ and }\mb e_1\prod_{i=1}^{n}A_n\mb e_1^t\sim c\varrho(A_1)\cdots\varrho(A_{n})\text{ as }n\rto.
   \end{align}
   Therefore, \eqref{epar} of Lemma \ref{smr} is proved.

  Finally, we turn to prove the first limit in \eqref{asal}. Notice that
  \begin{align}
    \sum_{k=1}^{n+1}&\mb e_1\prod_{i=k}^{n}A_i(1,\lambda_{n+1})^t
    =\sum_{k=1}^{n+1}\mb e_1\prod_{i=k}^{n}A_i\mb e_1^t+\lambda_{n+1}\sum_{k=1}^{n+1}\mb e_1\prod_{i=k}^{n}A_i\mb e_2^t\no\\
    &=\sum_{k=1}^{n+1}\mb e_1\prod_{i=k}^{n}A_i\mb e_1^t+\lambda_{n+1}b_n\sum_{k=1}^{n}\mb e_1\prod_{i=k}^{n-1}A_i\mb e_1^t\no.
  \end{align}
  Thus, we have
  \begin{align}\label{als}
    \frac{\sum_{k=1}^{n+1}\mb e_1\prod_{i=k}^{n}A_i(1,\lambda_{n+1})^t}{S_n}=\tilde Y_n+\lambda_{n+1}b_n\tilde Y_{n-1}\frac{S_{n-1}}{S_n}.
  \end{align}
Taking limits on both sides of \eqref{als} and using \eqref{tya} and \eqref{srl}, we obtain that
\begin{align}\label{asl}
    \lim_{n\rto}&\frac{\sum_{k=1}^{n+1}\mb e_1\prod_{i=k}^{n}A_i(1,\lambda_{n+1})^t}{S_n}\\
    &=\tilde\psi\z(1+(b-\theta)\z(1_{\{\varrho(A)<1\}}+\varrho(A)^{-1}1_{\{\varrho(A)\ge1\}}\y)\y).\no
  \end{align}
  What is left for us to show is the strict positivity of  the limit in \eqref{asl}. Since $\tilde\psi>0,$ it suffices to show
  \begin{align}\label{tps}
  \tau:=1+(b-\theta)\z(1_{\{\varrho(A)<1\}}+\varrho(A)^{-1}1_{\{\varrho(A)\ge1\}}\y)>0.
  \end{align}
  Clearly, if $b\ge \theta,$ we have $\tau>0.$ Thus we suppose next that $b<\theta.$
If we can show that
\begin{align*}
  b<\theta \text{ implies } \theta-b<\varrho(A),
\end{align*}
then $\tau>0.$ Indeed, if $b<\theta$, then $\theta-b<\varrho(A)$ if and only if
\begin{align}\label{lr}
  \theta-2b-a<\sqrt{(a+\theta)^2+4(b d-\theta a)}.
\end{align}
If $b<\theta\le a+2b,$ then clearly, \eqref{lr} is true.
If $\theta>a+2b,$ then \eqref{lr} holds if and only if $(\theta-2b-a)^2<(a+\theta)^2+4(b d-\theta a)$
or equivalently
$b-\theta+a<d,$ which holds trivially  since $b-\theta+a<0.$ Therefore, \eqref{tps} is proved and so is Lemma \ref{smr}. \qed

\subsection{Proof of Lemma \ref{bma}}
Suppose condition (B1) holds and $\tilde d_k\ge \varepsilon,\forall k\ge1$ for some $\varepsilon>0.$ Then we have $\sum_{k=2}^\infty|b_k-b_{k-1}|+|\theta_k-\theta_{k-1}|<\infty$ and $b_kd_k-a_k\theta_k>0, \forall k\ge1.$
 To begin with, we show that
\begin{align}\label{pma}
\varrho(M_1)\cdots (M_k)\sim c\varrho(A_1)\cdots \varrho(A_k) \text{ as }k\rto.
\end{align}
To this end, noticing that
 $\forall k\ge1,$ we have
  \begin{align*}
    \varrho(M_k)&=\frac{1}{2}\z(a_k+\theta_{k}+\sqrt{(a_k+\theta_{k})^2+4(b_kd_k-a_k\theta_k)}\y),\\
   \varrho(A_k) &=\frac{1}{2}\Big(a_k+\theta_k+b_k\Delta_k+\sqrt{\big(a_k+\theta_k+b_k\Delta_k\big)^2+4(b_kd_k-a_k\theta_k)}\Big),\no
  \end{align*}
  where $\Delta_k=:\frac{\theta_{k+1}}{b_{k+1}}-\frac{\theta_k}{b_k}.$
  It is easy to see that $\lim_{k\rto}\frac{\varrho(M_k)}{\varrho(A_k)}=1.$
 Also by some careful computation, taking condition (B1) into account, we get
  \begin{align*}
    \sum_{k=1}^\infty\varrho(M_k)^{-1}|\varrho(A_k)-\varrho(M_k)|\le c\sum_{k=2}^\infty|b_k-b_{k-1}|+|\theta_k-\theta_{k-1}|<\infty.
  \end{align*}
  Consequently, $
      \sum_{k=1}^\infty\log\z(1+\frac{\varrho(A_k)-\varrho(M_k)}{\varrho(M_k)}\y)
 $ is convergent. So we get \eqref{pma}.

 Next we show \begin{align}\label{eam}
   \sum_{k=1}^{n+1}\prod_{i=1}^{k-1} \varrho(A_i)^{-1}\sim c\sum_{k=1}^{n+1} \prod_{i=1}^{k-1} \varrho(M_i)^{-1}, \text{ as }n\rto.
 \end{align}
 Indeed, by \eqref{pma}, $\sum_{k=1}^{n+1}\prod_{i=1}^{k-1} \varrho(A_i)^{-1}$ and $\sum_{k=1}^{n+1} \prod_{i=1}^{k-1} \varrho(M_i)^{-1}$ converge or diverge simultaneously.
 If they are convergent,  then \eqref{eam} holds trivially. Otherwise, if they are divergent, then
$$\lim_{n\rto}\frac{\sum_{k=1}^{n+1} \prod_{i=1}^{k-1} \varrho(M_i)^{-1}}{\sum_{k=1}^{n+1}\prod_{i=1}^{k-1} \varrho(A_i)^{-1}}
=\lim_{n\rto}\frac{\varrho(M_1)^{-1}\cdots\varrho(M_{n})^{-1}}{\varrho(A_1)^{-1}\cdots\varrho(A_{n})^{-1}}=c$$
which implies \eqref{eam}.

As a result, we get Lemma \ref{bma} by taking \eqref{pma} and \eqref{eam} together.\qed

\section{Products of positive matrices and the tails of continued fractions}\label{pmp}
The main purpose of this section is to prove  Proposition \ref{tmc}. The proof relies heavily on various properties of the continued fractions and their tails. Keep also in mind that
  under condition (B1),  we have \begin{align*}
   &B_k:=\left(
         \begin{array}{cc}
           a_k & b_k \\
           d_k & 0 \\
         \end{array}
       \right)\rightarrow\left(
         \begin{array}{cc}
           a& b \\
           d & 0 \\
         \end{array}
       \right)=:B \text{ as }k\rto,\\
    &\varrho(B_k)=\frac{\sqrt{a_k^2+4b_kd_k}+a_k}{2}\rightarrow \frac{\sqrt{a^2+4bd}+a}{2}=:\varrho\equiv\varrho(B) \text{ as }k\rto.
 \end{align*}

\subsection{Matrix products, continued fractions and their approximants}
  \subsubsection{Notations of continued fractions}
 To begin with, we introduce some notations of  the continued fractions.
Let $\beta_k,\alpha_k>0,k\ge 1$ be certain  numbers. For $1\le k\le n,$ we denote by
\begin{equation}\label{aprx}
\xi_{k,n}\equiv\frac{\beta_k}{\alpha_k}\begin{array}{c}
                                \\
                               +
                             \end{array}\frac{\beta_{k+1}}{ \alpha_{k+1}}\begin{array}{c}
                                \\
                               +\cdots+
                             \end{array}\frac{\beta_n}{\alpha_n}:=\dfrac{\beta_k}{\alpha_k+\dfrac{\beta_{k+1}}{\alpha_{k+1}+_{\ddots_{\textstyle +\frac{\textstyle\beta_{n}}{\textstyle\alpha_{n} } }}}},
\end{equation}
 the $(n-k+1)$-th approximant of a continued fraction
 \begin{align}\label{xic}
   &\xi_k:=\frac{\beta_{k}}{\alpha_{k}}\begin{array}{c}
                                \\
                               +
                             \end{array}\frac{\beta_{k+1}}{\alpha_{k+1 }}\begin{array}{c}
                                \\
                               +
                             \end{array}\frac{\beta_{k+2}}{\alpha_{k+2}}\begin{array}{c}
                                \\
                               +\cdots
                             \end{array}.
\end{align}
If  $\lim_{n\rto}\xi_{k,n}$ exists, then we say that the continued fraction $\xi_k$  is convergent and its value is defined as $\lim_{n\rto}\xi_{k,n}.$ If \begin{align}
  \exists C>0 \text{ such that } \forall k\ge1,\ C^{-1}\le {\beta_k}/{\alpha_k}\le C,\label{ssc}
 \end{align}
 then  by Seidel-Stern Theorem (see Lorentzen and Waadeland \cite[Theorem 3.14]{lw}), for any $k\ge1,$  $\xi_k$ is convergent.
 In the literature, $\xi_k,k\ge1$ in \eqref{xic}  are usually called the tails of the continued fraction $\xi_1:={\frac{\beta_{1}}{\alpha_{1}}}_{+}\frac{\beta_{2}}{\alpha_{2 }}_{+\cdots}$
and $h_k:= \frac{\beta_k}{\alpha_{k-1}}_{+}\frac{\beta_{k-1}}{ \alpha_{k-2}}_{+\cdots+}\frac{\beta_2}{\alpha_1},k
                             \ge2$
 are  referred to as the critical tails of $\xi_1.$

\subsubsection{Product of matrices expressed in terms of the approximants of continued fractions}

Now, for $1\le k\le n,$ set
 \begin{align}\label{xiy}
   y_{k,n}:=\mathbf e_1 B_k\cdots B_{n}\mathbf e_1^t \text{ and }\xi_{k,n}:=\frac{y_{k+1,n}}{y_{k,n}}.
 \end{align} Noting that the empty product equals identity, thus $y_{n+1,n}=1.$ Therefore,
 \begin{align}
\xi_{k,n}^{-1}\cdots \xi_{n,n}^{-1}&=y_{k,n}=\mathbf e_1 B_k\cdots B_{n}\mathbf e_1^t,\label{cpn}\\
\label{axn}\sum_{k=1}^{n+1} \mathbf e_1B_k\cdots B_n\mathbf e_1^t&=\sum_{k=1}^{n+1} \xi_{k,n}^{-1}\cdots \xi_{n,n}^{-1}=\frac{\sum_{k=1}^{n+1} \xi_{1,n}\cdots \xi_{k-1,n}}{\xi_{1,n}\cdots \xi_{n,n}}.
 \end{align}
 \begin{lemma}\label{axc}For $1\le k\le n,$ $\xi_{k,n}$ defined in \eqref{xiy} coincides with the one in \eqref{aprx} with $\beta_k=b_k^{-1}d_{k+1}^{-1}$ and $\alpha_k=a_kb_k^{-1}d_{k+1}^{-1}.$
    \end{lemma}
    \proof  Clearly, $\xi_{n,n}=\frac{1}{y_{n,n}}=\frac{1}{a_n}=\frac{b_n^{-1}d_{n+1}^{-1}}{a_nb_n^{-1}d_{n+1}^{-1}}.$ For $1\le k< n,$ note that
\begin{align}\label{ix}
  \xi_{k,n}&=\frac{y_{k+1,n}}{y_{k,n}}=\frac{\mathbf e_1 B_{k+1}\cdots B_{n}\mathbf e_1^t}{\mathbf e_1 B_k\cdots B_{n}\mathbf e_1^t}=\frac{\mathbf e_1 B_{k+1}\cdots B_{n}\mathbf e_1^t}{(a_k\mathbf e_1+b_k\mathbf e_2) B_{k+1}\cdots B_{n}\mathbf e_1^t}\\
  &=\frac{1}{a_k+b_k\frac{\mathbf e_2 B_{k+1}\cdots B_{n}\mathbf e_1^t}{\mathbf e_1 B_{k+1}\cdots B_{n}\mathbf e_1^t}}=\frac{1}{a_k+b_kd_{k+1}\frac{\mathbf e_1 B_{k+2}\cdots B_{n}\mathbf e_1^t}{\mathbf e_1 B_{k+1}\cdots B_{n}\mathbf e_1^t}}\no\\
  &=\frac{b_k^{-1}d_{k+1}^{-1}}{a_kb_k^{-1}d_{k+1}^{-1}+\xi_{k+1,n}}.\no
  \end{align}
  Thus, the lemma can be proved by iterating \eqref{ix}. \qed

With
 \eqref{cpn} and \eqref{axn} in hands, Lemma \ref{xias} below, whose proof will be postponed to the end of this subsection,  is crucial to the proof of Proposition \ref{tmc}.

\begin{lemma}\label{xias} Let $\xi_k,k\ge1,$ $\xi_{k,n},1\le k\le n$ be the ones defined in \eqref{aprx} and \eqref{xic} with $\beta_k=b_k^{-1}d_{k+1}^{-1}$ and $\alpha_k=a_kb_k^{-1}d_{k+1}^{-1},k\ge1.$ Suppose condition (B1) and one of the conditions (B2)$_a,$ (B2)$_b$ and  (B2)$_c$ hold. Then as $n\rto$, we have
\begin{align}
 &\xi_{1,n}\cdots \xi_{n,n}\sim c\varrho(B_1)^{-1}\cdots\varrho(B_n)^{-1},  \label{xnr}\\
    &\xi_{1,n}\cdots \xi_{n,n}\sim c\xi_{1}\cdots \xi_{n}, \label{xnp}\\
      &\sum_{k=1}^{n+1}\xi_{1,n}\cdots\xi_{k-1,n}\sim c\sum_{k=1}^{n+1}\xi_{1}\cdots \xi_{k-1}. \label{sxnn}
 \end{align}
  \end{lemma}

\subsubsection{Lower and upper bounds for the product of approximants of continued fractions by the one of the tails}

The proof of Lemma \ref{xias} depends heavily on the following lemma which studies various properties of the continued fractions and their approximants.

 \begin{lemma}\label{xp} Let $\xi_{k,n}$ and $\xi_k$ be the ones in \eqref{aprx} and \eqref{xic}. Suppose that $\alpha_k,\beta_k>0, \forall k\ge1$ and \eqref{ssc} is satisfied.
 Then we have
 \begin{align}
 &\xi_{k,n}\rightarrow\xi_k\in(0,\infty), \text{ as } n\rto,\label{xil}\\
  & \xi_{k,n}\Big\{\begin{array}{cc}
     <\xi_k,& \text{if } n-k+1 \text{ is even,} \\
     >\xi_k,&\text{if } n-k+1 \text{ is odd,}
   \end{array},1\le k\le n,\label{pxi}\\
   &\xi_{k,n}\xi_{k+1,n}\Big\{\begin{array}{cc}
     >\xi_k\xi_{k+1},& \text{if } n-k+1 \text{ is even,} \\
     <\xi_k\xi_{k+1},&\text{if } n-k+1 \text{ is odd,}
   \end{array}1\le k\le n-1, \label{txi}\\
   &\xi_k\cdots\xi_n\le \xi_{k,n}\cdots\xi_{n,n} \le \xi_1\cdots\xi_{n-1}\beta_n/\alpha_n, n\ge k\ge1.\label{xud}
 \end{align}
  Furthermore if we assume in addition  $\alpha_k\rightarrow \alpha>0,$ $\beta_k\rightarrow\beta>0.$ Then  \begin{align}\label{ckx}
    &\lim_{k\rto}\xi_k=\xi:=\frac{\sqrt{\alpha^2+4\beta}-\alpha}{2}>0,\\
     &  \xi\le1\Rightarrow\sum_{k=1}^{n}\xi_{1,n}\cdots\xi_{k,n}\sim \sum_{k=1}^{n}\xi_1\cdots\xi_{k}.\label{sxud}
    \end{align}
  \end{lemma}
\proof Applying Seidel-Stern Theorem (see Lorentzen and Waadeland \cite[Theorem 3.14]{lw}), we get \eqref{xil} and with \eqref{xil} in hand,  \eqref{pxi} is a direct consequence of \cite[Thoerem 3.12]{lw}. Next, note that by \eqref{aprx} and \eqref{xic}, we have
\begin{align*}
  \xi_{k,n}\xi_{k+1,n}=\beta_k-\alpha_k\xi_{k,n} \text{ and } \xi_k\xi_{k+1}=\beta_k-\alpha_k\xi_k,
\end{align*}
respectively. Consequently, \begin{align*}
  \xi_{k,n}\xi_{k+1,n}- \xi_k\xi_{k+1}=\alpha_k(\xi_k-\xi_{k,n})
\end{align*}
which, together with \eqref{pxi}, implies  \eqref{txi}.

Next, we proceed to prove \eqref{xud}. We prove only the case $k$ is odd  and $n$ is  even, since the other cases  can be proved  similarly.  Noticing that $\xi_{n,n}=\frac{\beta_n}{\alpha_n},$
then it follows from \eqref{pxi} and \eqref{txi} that
\begin{align*}
  \xi_{k,n}\cdots\xi_{n,n}&=\xi_{k,n} \underline{\xi_{k+1,n}\xi_{k+2,n}}\cdots \underline{\xi_{n-2,n}\xi_{n-1,n}}\xi_{n,n}\le \xi_k\underline{\xi_{k+1}\xi_{k+2}}\cdots\underline{\xi_{n-2}\xi_{n-1}}\frac{\beta_n}{\alpha_n},\\
  \xi_{k,n}\cdots\xi_{n,n}&=\underline{\xi_{k,n}\xi_{k+1,n}}\cdots \underline{\xi_{n-1,n}\xi_{n,n}}\ge \underline{\xi_k\xi_{k+1}}\cdots\underline{\xi_{n-1}\xi_{n}}
\end{align*} which imply \eqref{xud}.

The convergence of $\xi_k\rightarrow\xi>0$ in \eqref{ckx} follows from the convergence of the tails of the limit period continued fraction, see for example, Lorentzen \cite[(4.2) on page 81]{lor}.

Finally, we assume $\xi\le1$  to show \eqref{sxud}.  We claim that $\exists N_0>0$ such that
\begin{align}
  n-k+1\text{ is even } &\Rightarrow \frac{\xi_{k,n}-\xi_k}{\xi_{k+1}-\xi_{k+1,n}}< r, \forall n>k>N_0, \label{udo}\\
  n-k+1\text{ is odd } &\Rightarrow \frac{\xi_{k,n}-\xi_k}{\xi_{k+2,n}-\xi_{k+2}}<r^2,\forall n>k>N_0\label{udt}
\end{align}
where $0<r<1$ is a proper number.

 In fact, since $\xi_k\rightarrow\xi\le1$ and $\alpha_k\rightarrow\alpha>0,$ then
 $\frac{\xi_{k}}{\alpha_k+\xi_{k+1}}\rightarrow\frac{\xi}{\alpha+\xi}<1$ as $k\rto.$
 As a result,  for some proper number $0<r<1,$ $\exists N_0>0$ such that
 $\frac{\xi_{k}}{\alpha_k+\xi_{k+1}}<r, \forall k\ge N_0.$
 On the other hand,
 it follows from \eqref{aprx} and \eqref{xic} that
 \begin{align*}
 \frac{\xi_{k,n}-\xi_k}{\xi_{k+1}-\xi_{k+1,n}}=\frac{\xi_{k,n}}{\alpha_k+\xi_{k+1}}=\frac{\xi_{k}}{\alpha_k+\xi_{k+1,n}}.
 \end{align*}
 If $n-k+1$ is even, then by \eqref{pxi}, $\xi_{k+1,n}>\xi_{k+1}.$ Thus
 \begin{align*}
 \frac{\xi_{k,n}-\xi_k}{\xi_{k+1}-\xi_{k+1,n}}=\frac{\xi_{k}}{\alpha_k+\xi_{k+1,n}}<\frac{\xi_{k}}{\alpha_k+\xi_{k+1}}<r, \ \forall n>k>N_0.
 \end{align*}
  If  $n-k+1$ is odd, then by \eqref{txi}, we have $\xi_{k,n}\xi_{k+1}<\xi_k\xi_{k+1}.$ Therefore,
 \begin{align*}
 \frac{\xi_{k,n}-\xi_k}{\xi_{k+2,n}-\xi_{k+2}}
 &=\frac{\xi_{k,n}}{\alpha_k+\xi_{k+1}}\frac{\xi_{k+1,n}}{\alpha_{k+1}+\xi_{k+2}}\\
 &<\frac{\xi_{k}}{\alpha_k+\xi_{k+1}}\frac{\xi_{k+1}}{\alpha_{k+1}+\xi_{k+2}}<r^2,\ \forall n\ge k>N_0.
 \end{align*}
 We thus finish  proving the claim.

Now, we begin to prove \eqref{sxud}. Again, we deal with the case $n$ is an even number only, since another case follows similarly. Assume that $n$ is an even number. Applying \eqref{pxi}, \eqref{txi} and checking carefully, we get
\begin{align}\label{ubd}
  \sum_{k=1}^{n}&\xi_{1,n}\cdots\xi_{k,n}=\xi_{1,n}+\xi_{1,n}\underline{\xi_{2,n}\xi_{3,n}}+...+\xi_{1,n} \underline{\xi_{2,n}\xi_{3,n}}\cdots \underline{\xi_{n-2,n}\xi_{n-1,n}}\\
&+\xi_{1,n}\xi_{2,n}+\xi_{1,n}\underline{\xi_{2,n}\xi_{3,n}}\xi_{4,n}+...+\xi_{1,n} \underline{\xi_{2,n}\xi_{3,n}}\cdots \underline{\xi_{n-2,n}\xi_{n-1,n}}\xi_{n,n}\no\\
&\le \sum_{k=1}^{n}\xi_1\cdots\xi_{k}+\sum_{k=1}^{n/2}\prod_{i=1}^{2k-1}\xi_i(\xi_{2k,n}-\xi_{2k}).\no
\end{align}
Since $n-2k+1$ is odd, we have $\xi_{2k,n}-\xi_{2k}>0$ and thus the second summation in the rightmost hand of the above inequality is positive. Next we show that
\begin{align}\label{dps}
  \lim_{n\rto}\frac{\sum_{k=1}^{n/2}\prod_{i=1}^{2k-1}\xi_i(\xi_{2k,n}-\xi_{2k})}{\sum_{k=1}^{n}\xi_1\cdots\xi_{k}}=0.
\end{align}

To this end, fix $\varepsilon>0$ and let $N_1>0$ be an even number such that $r^{N_1}<\varepsilon.$
For convenience, we may assume $N_0$ in \eqref{udo} and \eqref{udt} is odd.
Then
\begin{align}\label{tsp}
  \sum_{k=1}^{n/2}\prod_{k=1}^{2k-1}\xi_i(\xi_{2k,n}-\xi_{2k})&=\sum_{k=1}^{\frac{N_0-1}{2}}
  +\sum_{k=\frac{N_0-1}{2}+1}^{\frac{n-N_1}{2}}+\sum_{k=\frac{n-N_1}{2}+1}^{\frac{n}{2}}\prod_{i=1}^{2k-1}\xi_i(\xi_{2k,n}-\xi_{2k})\\
  &=:\mathrm{(I) +(II)+(III).}\no
\end{align}
Since $\lim_{n\rto}\xi_{k,n}=\xi_k$ by \eqref{xil}, the first term (I) on the rightmost hand of \eqref{tsp} vanishes  as $n\rto.$
Using the fact $\lim_{k\rto}\xi_k \le 1$ and applying Lemma \ref{sig}, for each $(n-N_1)\le i\le n/2,$ we have
$\frac{\xi_1\cdots\xi_{n-i}}{\sum_{k=1}^{n}\xi_1\cdots\xi_{k}}\rightarrow 0$ as $n\rto.$ Note also that $|\xi_{k,n}-\xi_{k}|<2\beta_k/\alpha_k<C$ for some universal constant $C>0.$
Therefore, we conclude that
$$\lim_{n\rto}\mathrm{(III)}/\sum_{k=1}^{n}\xi_1\cdots\xi_{k}=0.$$
Now we consider the term (II) on the rightmost hand of \eqref{tsp}. It follows from \eqref{udt} that $\xi_{2k,n}-\xi_{2k}<r^{n-2k}(\xi_{n,n}-\xi_n)\le Cr^{n-2k}$ for some constant $C>0$ independent of $n$ and $k.$
Then we have
\begin{align*}
\frac{\mathrm{(II)}}{\sum_{k=1}^{n}\xi_1\cdots\xi_{k}}&= \sum_{k=(N_0-1)/2+1}^{(n-N_1)/2}\prod_{i=1}^{2k-1}\xi_i(\xi_{2k,n}-\xi_{2k})\Big/ \sum_{k=1}^{n}\xi_1\cdots\xi_{k}\\
  &\le \sum_{k=(N_0-1)/2+1}^{(n-N_1)/2} Cr^{n-2k} <\sum_{k=N_1/2}^{\infty} Cr^{2k}\\
  &=Cr^{N_1}/(1-r^2)<C\varepsilon/(1-r^2).
\end{align*}
Since $\varepsilon$ is arbitrary, we get $\lim_{n\rto}\frac{\mathrm{(II)}}{\sum_{k=1}^{n}\xi_1\cdots\xi_{k}}=0.$ Therefore,
we come to the conclusion that \eqref{dps} is true. As a consequence, dividing by $\sum_{k=1}^{n}\xi_1\cdots\xi_{k}$ on both sides of \eqref{ubd} and taking the upper limit, we conclude that
\begin{align}\label{upb}
  \varlimsup_{n\rto}\frac{\sum_{k=1}^{n}\xi_{1,n}\cdots\xi_{k,n}}{\sum_{k=1}^{n}\xi_1\cdots\xi_{k}}=1.
\end{align}
For a lower limit, from \eqref{pxi}, \eqref{txi} and \eqref{udo} we get
\begin{align*}
  \sum_{k=1}^{n}&\xi_{1,n}\cdots\xi_{k,n}=\sum_{k=1}^{n/2}\underline{\xi_{1,n}\xi_{2,n}}\cdots \underline{\xi_{2k-3,n}\xi_{2k-2,n}}(\xi_{2k-1,n}-\xi_{2k-1})\\
 &+\sum_{k=1}^{n/2}\underline{\xi_{1,n}\xi_{2,n}}\cdots \underline{\xi_{2k-3,n}\xi_{2k-2,n}}\xi_{2k-1}+\sum_{k=1}^{n/2}\underline{\xi_{1,n}\xi_{2,n}}\cdots \underline{\xi_{2k-1,n}\xi_{2k,n}}\\
 &\ge \sum_{k=1}^n \xi_1\xi_2\cdots\xi_k+\bigg(\sum_{k=1}^{N_0/2}+\sum_{k=N_0/2+1}^{n/2}\underline{\xi_{1}\xi_{2}}\cdots \underline{\xi_{2k-3}\xi_{2k-2}}(\xi_{2k,n}-\xi_{2k})\bigg).
 \end{align*}
  Using \eqref{udt}, similar to \eqref{dps}, we have
  \begin{align*}\lim_{n\rto}\frac{\sum_{k=1}^{N_0/2}+\sum_{k=N_0/2+1}^{n/2}{\xi_{1}\xi_{2}}\cdots {\xi_{2k-3}\xi_{2k-2}}(\xi_{2k,n}-\xi_{2k})}{\sum_{k=1}^n \xi_1\xi_2\cdots\xi_k}=0.
      \end{align*}
  As a consequence, \begin{align}\label{lowb}
  \varliminf_{n\rto}\frac{\sum_{k=1}^{n}\xi_{1,n}\cdots\xi_{k,n}}{\sum_{k=1}^{n}\xi_1\cdots\xi_{k}}=1.
\end{align}
  Taking \eqref{upb} and \eqref{lowb} together, we get \eqref{sxud}.  The lemma is proved.
  \qed

 \subsubsection{Fluctuations of tail and critical tail of continued fractions}
For $k\ge1,$ let \begin{align*}
  f_k&=\frac{ b_k d_k^{-1}}{ a_kd_k^{-1}}\begin{array}{c}
                                \\
                               +
                             \end{array}\frac{ b_{k-1} d_{k-1}^{-1}}{ a_{k-1}d_{k-1}^{-1}}\begin{array}{c}
                                \\
                               +\cdots+
                             \end{array}\frac{ b_1 d_1^{-1}}{ a_1 d_1^{-1}},\\
\xi_k&=\frac{ b_k^{-1} d_{k+1}^{-1}}{ a_k b_k^{-1} d_{k+1}^{-1}}\begin{array}{c}
                                \\
                               +
                             \end{array}\frac{ b_{k+1}^{-1} d_{k+2}^{-1}}{ a_{k+1}b_{k+1}^{-1}d_{k+2}^{-1}}\begin{array}{c}
                                \\
                               +\cdots.
                             \end{array}
\end{align*}
Set also \begin{align*}
  \varepsilon^f_k&=f_k-b_{k+1}\varrho(B_{k+1})^{-1},\ \varepsilon^\xi_k=\xi_k-\varrho(B_{k})^{-1}, k\ge1,\\
  \delta_k^f&=b_kd_k^{-1}-b_{k+1}\varrho(B_{k+1})^{-1}(a_kd_k^{-1}+b_{k}\varrho(B_{k})^{-1}),k\ge2,\\
  \delta_k^\xi&=b_k^{-1}d_{k+1}^{-1}-\varrho(B_k)^{-1}(a_kb_k^{-1}d_{k+1}^{-1}+\varrho(B_{k+1})^{-1}),k\ge1.
\end{align*}
Suppose that condition (B1) holds. Then
by the theory of  convergence of limit periodic continued fractions (see  Lorentzen \cite[(4.2) on page 81]{lor} and  Lorentzen and Waadeland \cite[Theorem 4.13, page 188]{lw}), we have
\begin{align}\label{xkl}
  f_k\rightarrow b\varrho^{-1},\ \xi_k\rightarrow\varrho^{-1} \text{ and hence }\varepsilon_k^f \rightarrow0,\ \varepsilon_k^\xi \rightarrow0 \text{ as }k\rto.
\end{align}
Lemma \ref{dxf} below studies the fluctuations of $\varepsilon_k^f$ and $\varepsilon_k^\xi$ $k\ge1.$
\begin{lemma}\label{dxf}Suppose that condition (B1) and one of the conditions (B2)$_{a},$ (B2)$_{ b}$  and  (B2)$_{c}$ hold. Then there exists some number $q$ with $|q|\le1$ such that
\begin{align}\label{qdd}
  &\lim_{k\rto}{\delta_{k+1}^f}/{\delta_k^f}=\lim_{k\rto}{\delta_{k+1}^\xi}/{\delta_k^\xi}=q,\\
 \label{qff} &\lim_{k\rto}\frac{\varepsilon_{k+1}^\xi}{\varepsilon_{k}^\xi}=q,\  \lim_{k\rto}\frac{\varepsilon_{k+1}^f}{\varepsilon_{k}^f}=q\text{ or }-\frac{b\varrho^{-1}}{1+b\varrho^{-1}}.
\end{align}
\end{lemma}
\proof If \eqref{qdd} holds, then \eqref{qff} is a direct consequence of \cite[Lemma 4]{hs20} and Lorentzen \cite[Theorem 6.1]{lor}. Therefore, it is sufficient to prove \eqref{qdd}.

Suppose condition (B1) holds. Then by \eqref{xkl}, $\lim_{k\rto} \delta_k^\xi=\lim_{k\rto} \delta_k^f=0.$ Thus, if the limits in \eqref{qdd} exist, we must have $|q|\le1.$ It remains to show that the limits in \eqref{qdd} exist and are equal.
To this end, notice that by some direct computation, we get
\begin{align*}
  \delta_k^f=\frac{b_kd_k^{-1}\Delta_k^f}{\frac{a_{k+1}}{b_{k+1}}+\sqrt{\z(\frac{a_{k+1}}{b_{k+1}}\y)^2+4\frac{d_{k+1}}{b_{k+1}}}}\text{ and }
  \delta_k^\xi=\frac{b_k^{-1}d_{k+1}^{-1}\Delta_k^\xi}{\frac{a_{k}}{b_{k}}+\sqrt{\z(\frac{a_{k}}{b_{k}}\y)^2+4\frac{d_{k}}{b_{k}}}}
\end{align*}
where
$$\Delta_k^f=\frac{a_{k+1}}{b_{k+1}}-\frac{a_{k}}{b_{k}}+\frac{\z(\frac{a_{k+1}}{b_{k+1}}-\frac{a_{k}}{b_{k}}\y)\z(\frac{a_{k+1}}{b_{k+1}}+\frac{a_{k}}{b_{k}}\y)+4\z(\frac{d_{k+1}}{b_{k+1}}-\frac{d_{k}}{b_{k}}\y)}{\sqrt{\z(\frac{a_{k+1}}{b_{k+1}}\y)^2+4\frac{d_{k+1}}{b_{k+1}}}+\sqrt{\z(\frac{a_{k}}{b_{k}}\y)^2+4\frac{d_{k}}{b_{k}}}}$$
and $$\Delta_k^\xi=\frac{a_{k+1}}{b_{k+1}}-\frac{a_{k}}{b_{k}}-\frac{\z(\frac{a_{k+1}}{b_{k+1}}-\frac{a_{k}}{b_{k}}\y)\z(\frac{a_{k+1}}{b_{k+1}}+\frac{a_{k}}{b_{k}}\y)+4\z(\frac{d_{k+1}}{b_{k+1}}-\frac{d_{k}}{b_{k}}\y)}{\sqrt{\z(\frac{a_{k+1}}{b_{k+1}}\y)^2+4\frac{d_{k+1}}{b_{k+1}}}+\sqrt{\z(\frac{a_{k}}{b_{k}}\y)^2+4\frac{d_{k}}{b_{k}}}}.$$
Therefore, by condition (B1) we have
$$\lim_{k\rto}{\delta_{k+1}^f}/{\delta_k^f}=\lim_{k\rto}{\Delta_{k+1}^f}/{\Delta_k^f},\lim_{k\rto}{\delta_{k+1}^\xi}/{\delta_k^\xi}=\lim_{k\rto}{\Delta_{k+1}^\xi}/{\Delta_k^\xi}. $$
Suppose now condition (B2)$_{a}$ holds. Then the limits
$$\lim_{k\rto}{\delta_{k+1}^f}/{\delta_k^f}=\lim_{k\rto}{\delta_{k+1}^\xi}/{\delta_k^\xi}=\lim_{k\rto}\frac{d_{k+2}/b_{k+2}-d_{k+1}/b_{k+1}}{d_{k+1}/b_{k+1}-d_{k}/b_{k}}$$
exist. Next, suppose that condition (B2)$_{b}$ holds. Then the limits
$$\lim_{k\rto}{\delta_{k+1}^f}/{\delta_k^f}=\lim_{k\rto}{\delta_{k+1}^\xi}/{\delta_k^\xi}=\lim_{k\rto}\frac{a_{k+2}/b_{k+2}-a_{k+1}/b_{k+1}}{a_{k+1}/b_{k+1}-a_{k}/b_{k}}$$ exist. Finally, suppose condition  (B2)$_{ c}$ holds. If $$\tau:=\lim_{k\rto}\frac{d_{k+1}/b_{k+1}-d_{k}/b_k}{a_{k+1}/b_{k+1}-a_{k}/b_k}\ne \frac{-a\pm \sqrt{a^2+4bd}}{2b}  $$ is finite, then
\begin{align*}
 &\lim_{k\rto}\frac{\Delta_k^f}{\frac{a_{k+1}}{b_{k+1}}-\frac{a_{k}}{b_{k}}} =1-\frac{a/b+2\tau}{\sqrt{(a/b)^2+4d/b}}\ne 0,\\
 &\lim_{k\rto}\frac{\Delta_k^\xi}{\frac{a_{k+1}}{b_{k+1}}-\frac{a_{k}}{b_{k}}} =1+\frac{a/b+2\tau}{\sqrt{(a/b)^2+4d/b}}\ne 0
\end{align*} and consequently, the limits
\begin{align*}
  \lim_{k\rto}{\delta_{k+1}^f}/{\delta_k^f}=\lim_{k\rto}{\delta_{k+1}^\xi}/{\delta_k^\xi}=\lim_{k\rto}\frac{a_{k+2}/b_{k+2}-a_{k+1}/b_{k+1}}{a_{k+1}/b_{k+1}-a_{k}/b_{k}}
\end{align*} exist. Otherwise, if $\tau=\infty,$ then $\lim_{k\rto}\frac{a_{k+1}/b_{k+1}-a_{k}/b_k}{d_{k+1}/b_{k+1}-d_{k}/b_k}=0$ and hence
$$\lim_{k\rto}{\delta_{k+1}^f}/{\delta_k^f}=\lim_{k\rto}{\delta_{k+1}^\xi}/{\delta_k^\xi}=\lim_{k\rto}\frac{d_{k+2}/b_{k+2}-d_{k+1}/b_{k+1}}{d_{k+1}/b_{k+1}-d_{k}/b_{k}}$$
exist. The lemma is proved.\qed


We are now ready to prove Lemma \ref{xias}.

\subsubsection{ Proof of Lemma \ref{xias}}  Suppose now condition (B1) and one of the conditions (B2)$_a,$ (B2)$_b$ and  (B2)$_c$ hold. Then
    Lemma \ref{dxf} ensures us to apply \cite[Theorem 1]{hs20} to $B_k^t,k\ge1$ to yield that for $i,j\in\{1,2\}$ and $m\ge 1,$ there exists a number $0<c(m)<\infty$ such that
\begin{align}\label{xer}
 \mb e_iB_m\cdots B_n\mb e_j^t =\mb e_jB_n^t\cdots B_m^t\mb e_i^t\sim c(m)\varrho(B_m)\cdots \varrho(B_n) \text{ as }n\rto.
\end{align}

Here, we remark that when considering $B_k^t,k\ge1,$ the conditions  (B2)$_a,$ (B2)$_b$ and  (B2)$_c$ in this paper are slightly different from the counterpart in \cite{hs20}.  But  the key role which (B2)$_a,$ (B2)$_b$ and  (B2)$_c$ in \cite{hs20} play is to show that $\lim_{k\rto}\frac{\delta_{k+1}}{\delta_k}$ exists. Lemma \ref{dxf} showed that  $\lim_{k\rto}\frac{\delta_{k+1}}{\delta_k}$ does exist under one of the conditions (B2)$_a,$ (B2)$_b$ and  (B2)$_c$ of this paper.  So, we can apply \cite[Theorem 1]{hs20} to $B_k^t,k\ge1$ to get \eqref{xer}.

 By Lemma \ref{axc}, for $1\le k\le n,$ $\xi_{k,n}$ in \eqref{aprx} coincides with the one in \eqref{xiy}, if $\beta_k=b_k^{-1}d_{k+1}^{-1}$ and $\alpha_k=a_kb_k^{-1}d_{k+1}^{-1}.$ Therefore from \eqref{cpn} and \eqref{xer} we get
\begin{align}\label{xerr}
  \xi_{1,n}\cdots \xi_{n,n}=\frac{1}{\mb e_1B_1\cdots B_n\mb e_1^t}\sim c\varrho(B_1)^{-1}\cdots \varrho(B_n)^{-1}, \text{ as }n\rto,
\end{align}
which finishes the proof of \eqref{xnr}.

  For the proof of \eqref{xnp}, noticing that
   $\lim_{k\rto}\beta_k=b^{-1}d^{-1}$ and $\lim_{k\rto}\alpha_k=ab^{-1}d^{-1},$ then \eqref{ssc} is fulfilled. Therefore, by Lemma \ref{xp} we have  for $k\ge1,$
  \begin{align*}
 &\xi_{k,n}\rightarrow\xi_k\in(0,\infty), \text{ as } n\rto,
   \end{align*}
      and for $n\ge1,$
    \begin{align}\label{xud2}
     &1\le \frac{\xi_{1,n}\cdots\xi_{n,n}}{\xi_1\cdots\xi_n} \le \frac{\beta_n}{\xi_n\alpha_n}.
    \end{align}
%
%
%
    To  prove \eqref{xnp}, in view of \eqref{xnr}, it suffices to show that
\begin{align}
 \xi_1\cdots\xi_n\sim c\varrho(B_1)^{-1}\cdots \varrho(B_n)^{-1} \text{ as }n\rto.\label{xirr}
\end{align}
To this end, write $x_n=\frac{\xi_1\cdots\xi_n}{\varrho(B_1)^{-1}\cdots\varrho(B_n)^{-1}},n\ge1.$ Noting that by  \eqref{xerr} and \eqref{xud2}, there exist some numbers $0< c_3\le c_4<\infty$ such that
\begin{equation}\label{lux}
  c_3\le x_k\le c_4, \forall k\ge1.
\end{equation}
 Notice that by Lemma \ref{dxf},
\begin{align}\label{xlq}
  \lim_{k\rto}\frac{\xi_{k+1}-\varrho(B_{k+1})^{-1}}{\xi_{k}-\varrho(B_{k})^{-1}}=q,
\end{align}
for some number $q$ with $|q|\le1.$

At first, we suppose  $|q|<1.$  Then it follows from \eqref{xlq} that $\sum_{k=1}^{\infty}|\xi_k-\varrho(B_k)^{-1}|<\infty.$ Since $\lim_{n\rto}\varrho(B_n)=\varrho,$ taking \eqref{lux} into account, we have
\begin{align}\label{asc}
  \sum_{n=1}^\infty|x_{n+1}-x_{n}|&=\sum_{n=1}^{\infty}\varrho(B_{n+1})x_n|\xi_{n+1}-\varrho(B_{n+1})^{-1}|\\
  &\le c\sum_{n=1}^{\infty}|\xi_{n+1}-\varrho(B_{n+1})^{-1}|<\infty.\no
\end{align}
Therefore, by \eqref{lux} and \eqref{asc}, we have $\lim_{k\rto}x_k=c$ for some $c>0.$

Secondly, we suppose $q=1.$ Then, there exists a number $k_3>0$ such that
\begin{align}\label{mo}
  \text{ either }\xi_k-\varrho(B_k)^{-1} \ge 0, \forall k\ge k_3 \text{ or } \xi_k-\varrho(B_k)^{-1} \le 0, \forall k\ge k_3.
\end{align}
But for $k\ge1$ we have
$ \frac{x_{k+1}}{x_k}=\frac{\xi_{k+1}}{\varrho(B_{k+1})^{-1}},$
so that by \eqref{mo}, we must have
  either  $x_{k+1}\le x_k, \forall k\ge k_3$  or $x_{k+1}\ge x_k, \forall k\ge k_3.$ That is, $x_k,k\ge k_3$ is monotone in $k.$ This fact together with \eqref{lux} implies $\lim_{k\rto}x_k=c$ for some number $c>0.$

  Finally, we suppose $q=-1.$ Then there exists a number $k_4>0$ such that $$\frac{\xi_{k+1}-\varrho(B_{k+1})^{-1}}{\xi_{k}-\varrho(B_{k})^{-1}}<0, \forall k\ge k_4.$$
   Thus, $\xi_{k}-\varrho(B_{k})^{-1},k\ge k_4$ converges to $0$ in an alternate manner as $k\rto.$ Therefore $\lim_{k\rto}x_k=\lim_{k\rto}\z(x_1+\sum_{i=2}^{k}(x_k-x_{k-1})\y)=c$ for some  $c>0.$

   We have shown that in any case,  $\lim_{k\rto}x_k=c$ for some  $c>0.$ Consequently,  \eqref{xirr} is proved and so is \eqref{xnp}.

   Finally, we give the proof of \eqref{sxnn}.
   If $\varrho\ge1,$ then by \eqref{xkl}, $\xi=\varrho^{-1}\le1.$  Thus applying Lemma \ref{xp}, from \eqref{sxud}, we get \eqref{sxnn}.

   Therefore, we suppose next $\varrho<1.$
   Since $\lim_{k\rto}\xi_k=\varrho^{-1}$ by \eqref{xkl}, applying Lemma
    \ref{sig} to $\xi_k,k\ge1,$ we get
    \begin{align}\label{xrl}
    \lim_{n\rto}\frac{\xi_1\cdots\xi_{n+1}}{\sum_{k=1}^{n+1}\xi_1\cdots\xi_{k-1}}
    =\z\{\begin{array}{cc}
      0, & \text{if }\varrho\ge 1, \\
    \varrho^{-1}-1, & \text{if }\varrho< 1.
    \end{array}
    \y.
    \end{align}
      It follows from \eqref{axn},  \eqref{xnp} and \eqref{xrl} that
   \begin{align*}
     \sum_{k=1}^{n+1} \mathbf e_1B_k\cdots B_n\mathbf e_1^t=\frac{\sum_{k=1}^{n+1} \xi_{1,n}\cdots \xi_{k-1,n}}{\xi_{1,n}\cdots \xi_{n,n}}\sim c\frac{\sum_{k=1}^{n+1}\xi_{1,n}\cdots\xi_{k-1,n}}{\sum_{k=1}^{n+1}\xi_1\cdots\xi_{k-1}},
   \end{align*}
   as $n\rto.$ Thus, to prove \eqref{sxnn}, it is sufficient to show that
   \begin{align}\label{lac}
     \lim_{n\rto}\sum_{k=1}^{n+1} \mathbf e_1B_k\cdots B_n\mathbf e_1^t=c
   \end{align}
   for some number $c>0.$
   Notice that by \eqref{cpn} and \eqref{xud}, we have for $n\ge m\ge1$
   \begin{align}\label{lumn}
 c_5\xi_m^{-1}\cdots \xi_n^{-1} <\mb e_1B_m\cdots B_n\mb e_1^t <c_6\xi_m^{-1}\cdots \xi_n^{-1}
\end{align}
   where $0<c_5<c_6<\infty$ are some numbers independent of $n$ and $m.$
   Since $\lim_{k\rto}\xi_k^{-1}=\varrho<1,$ fixing   $0\le \varepsilon_0 <(1-\varrho)\wedge\varrho,$ there exists a number $n_0>0$ such that
\begin{align}
  \varrho-\varepsilon_0<\xi_k^{-1}<\varrho+\varepsilon_0, \forall k\ge n_0.\label{rep}
\end{align}
Now fix $\varepsilon>0$ and let $n_1$ be a number large enough so that $(\varrho+\varepsilon_0)^{n_1}<\varepsilon.$

For $n\ge1,$ use the notation $y_{k,n}=\mathbf e_1B_k\cdots B_n\mathbf e_1^t$ and write $\Gamma_n=\sum_{k=1}^{n+1} y_{k,n}.$
 We have
\begin{align}\label{gas}
  |\Gamma_{n+1}&-\Gamma_n|=\Big|y_{1,n+1}+\sum_{k=1}^{n}(y_{k+1,n+1}-y_{k,n})\Big|\\
  &\le y_{1,n}+\Big|\sum_{k=1}^{n_0-1}(y_{k+1,n+1}-y_{k,n})\Big|+\Big|\sum_{k=n_0}^{n-n_1}(y_{k+1,n+1}-y_{k,n})\Big| \no\\
  &\quad\quad\quad\quad\quad+\Big|\sum_{k=n-n_1+1}^{n}(y_{k+1,n+1}-y_{k,n})\Big|\no\\
  &=:\mathrm{(I)+(II)+(III)+(IV)}.\no
\end{align}
Since $\varrho<1,$ for each $k\ge1,$ by \eqref{xer} we have $\lim_{n\rto}y_{k,n}=0.$ Also, for any $m\ge 1,$ the fact $\lim_{n\rto}B_n=B$ implies that $y_{n-m,n}-y_{n-m+1,n+1}\rightarrow 0$ as $n\rto.$ Thus,  (I) (II) and (IV) on the rightmost hand of \eqref{gas} vanish as $n\rto.$ We claim that (III) on the rightmost of \eqref{gas} also vanishes. Indeed, from  \eqref{lumn} we obtain
$y_{k,n}\le c_6\xi_k^{-1}\cdots \xi_n^{-1}$ for all $n\ge k\ge1.$ Therefore, taking \eqref{rep} into account,
\begin{align}
  \mathrm{(III)}&\le \sum_{k=n_0}^{n-n_1}(|y_{k+1,n+1}|+|y_{k,n}|)\le 2 c_6 \sum_{k=n_0}^{n-n_1}(\varrho+\varepsilon_0)^{n-k+1}\no\\
  &\no\le 2c_6\sum_{k=n_1+1}^{\infty}(\varrho+\varepsilon_0)^k=2c_6\frac{(\varrho+\varepsilon_0)^{n_1+1}}{1-(\varrho+\varepsilon_0)}
  \le \frac{2c_6\varepsilon}{1-(\varrho+\varepsilon_0)}.
\end{align}
Thus we can conclude that
\begin{align}\label{glim}
  \lim_{n\rto}(\Gamma_{n+1}-\Gamma_n)=0.
\end{align}
But by \eqref{lumn},  we have
\begin{align}\label{galu}
  c_7<\Gamma_{n}<c_8, \forall n\ge 1,
\end{align}
where $0<c_7<c_8<\infty$ are some constants independent of $n.$
Therefore, taking \eqref{glim} and \eqref{galu} together, we get
\begin{align}\label{gml}
  \lim_{n\rto}\frac{\Gamma_{n+1}}{\Gamma_n}=1.
\end{align}
For $n\ge1,$ set  \begin{align*}
  f_n\equiv \frac{\mb e_1\prod_{k=1}^{n}B_{k}\mathbf e_2^t}{\mb e_1\prod_{k=1}^{n} B_{k}\mathbf e_1^t} \text{ and }H_n\equiv\sum_{k=1}^{n}\mb e_1\prod_{i=k}^{n} B_{i}(f_n\mathbf e_1^t-\mb e_2^t).
\end{align*}Since $\varrho<1,$ we can apply Lemma \ref{hs} (with $\theta=0$) to get
  \begin{align}\label{fhl}
    \lim_{n\rto}f_n=-\frac{\varrho_1(A)}{d} \text{ and }  \lim_{n\rto}H_n=-\frac{1}{d }\frac{\varrho_1(A)^2}{1-\varrho_1(A)}.
  \end{align}
  Noticing that
  \begin{align}
  H_n&=(\Gamma_n-1)f_n-\sum_{k=1}^n\mb e_1 B_k\cdots B_{n}\mb e_2^t \no\\
   &=(\Gamma_n-1)f_n-b_n\Big(\sum_{k=1}^{n}\mb e_1 B_k\cdots B_{n-1}\mb e_1^t\Big)=\Gamma_nf_n-b_n\Gamma_{n-1}-f_n, \no
  \end{align}
thus it follows from \eqref{gml} and \eqref{fhl} that
\begin{align*}
  \lim_{n\rto}\Gamma_n=\lim_{n\rto} \frac{H_n+f_n}{f_n-b_n\Gamma_{n-1}/\Gamma_n}=\frac{\varrho_1(B)}{(\varrho_1(B)+bd)(1-\varrho_1(B))}>0,
\end{align*}
where to show positivity of the limit we use the fact $\varrho<1.$  Consequently,  \eqref{lac} is proved. Thus \eqref{sxnn} is true and we complete the proof of Lemma \ref{xias}.
  \qed

\subsection{Proof of Proposition \ref{tmc}}
  We  are now ready to complete the proof of Proposition \ref{tmc}. Suppose condition (B1) and one of the conditions (B2)$_a,$ (B2)$_b$ and  (B2)$_c$ hold. Then \eqref{rac} is a direct consequence of \eqref{xer} above.  Furthermore,
      it follows from \eqref{axn} that
   \begin{align*}
    Y_n&=\frac{\varrho(B_1)^{-1}\cdots\varrho(B_n)^{-1}}{\xi_{1,n}\cdots\xi_{n,n}}\frac{\sum_{k=1}^{n+1}\xi_{1,n}\cdots \xi_{k-1,n}}{\sum_{k=1}^{n+1}\xi_{1}\cdots \xi_{k-1}} \frac{\sum_{k=1}^{n+1}\xi_{1}\cdots \xi_{k-1}}{\sum_{k=1}^{n+1}\varrho(B_1)^{-1}\cdots\varrho(B_{k-1})^{-1}}\\
     &=:T_1(n)\times T_2(n)\times T_3(n).
   \end{align*}
   Noticing that by \eqref{xnr}, \eqref{xnp} and \eqref{sxnn}, for $i=1,2,3,$ we have
\begin{align}\label{ttt}\lim_{n\rto}T_i(n)=T_i,\end{align}
   for some number $T_i\in (0,\infty).$ Therefore $\lim_{n\rto}Y_n=T_1T_2T_3.$
Consequently, with \begin{align}\label{psi}
  \psi:=T_1T_2T_3\in(0,\infty)
\end{align} we get
\eqref{sax}.  Proposition \ref{tmc} is proved. \qed

\section*{Appendix: Proof of Lemma \ref{egc} }

Since $\lim_{n\rto}\frac{r_n-r_{n+1}}{r_n^2}=c$ for some number $0<c<\infty,$ we must have $\lim_{n\rto}\frac{r_{n+1}}{r_n}=1$ and there exists a number $N_0$ such that $\forall n\ge N_0, r_{n+1}<r_n.$
Consequently, the condition (B1) holds.  Next we show that one of the conditions (B2)$_{ a},$ (B2)$_{ b}$ and (B2)$_{ c}$ holds.

{\it Case 1:  Assume that $a=b=\theta\ne d.$ } Recall that by assumption $b>0.$ Thus for all $k\ge1,$ $\frac{\tilde a_k}{\tilde b_k}=\frac{a_k}{b_k}+\frac{\theta_{k+1}}{b_{k+1}}=\frac{a+r_k}{b+r_k}+\frac{\theta+r_{k+1}}{b+r_{k+1}}\equiv 2.$ That is, $\forall k\ge1, \frac{\tilde a_k}{\tilde b_k}=\frac{\tilde a_{k+1}}{\tilde b_{k+1}}.$ Moreover, notice that
\begin{align}\label{fdb}
  \frac{\tilde d_k}{\tilde b_k}&=\frac{d_kb_k-a_k\theta_k}{b_k^2}=\frac{(d+r_k)(b+r_k)-(a+r_k)(\theta+r_k)}{(b+r_k)^2}\\
  &=\frac{bd-a\theta+(b+d-a-\theta)r_k}{(b+r_k)^2}, k\ge1.\no
\end{align}
Writing $$g(x):=\frac{bd-a\theta+(b+d-a-\theta)x}{(b+x)^2},x\ge 0,$$ we have
\begin{align}\label{gd}
  g'(x)=\frac{b(b+d-a-\theta)+2(a\theta-bd)-(b+d-a-\theta)x}{(b+x)^3}.
\end{align}
Since $a=b=\theta\ne d,$ we get $g'(0)=\frac{b(b-d)}{b^3}\ne 0.$  Thus for some $k_0>0,$ $\frac{\tilde d_k}{\tilde b_k},k\ge k_0$ is monotone in $k.$ That is $\tilde d_{k}/\tilde b_{k}\ne\tilde d_{k+1}/\tilde b_{k+1},k\ge k_0.$

One the other hand,
for any numbers $\alpha,\beta\in \mathbb R$  with $\beta=0,$ we have
\begin{align}\label{asx}
  \frac{\alpha+x}{\beta+x}=\frac{\alpha}{\beta}+\frac{\beta-\alpha}{\beta^2}x-\frac{\beta-\alpha}{\beta^3}x^2+o(x^2),
\end{align}
as $x\rightarrow0.$ Thus, we get
\begin{align*}
  \tilde d_k/\tilde b_k=\frac{d+r_k}{b+r_k}-1=\frac{d}{b}-1 +\frac{b-d}{b^2}r_k-\frac{b-d}{b^3}r_k^2+o(r_k^2).
\end{align*}
As a result, using the fact $r_{n}\sim r_{n+1}$ and $r_n-r_{n+1}\sim cr_n^2$ as $n\rto,$ we get
$$\lim_{k\rto}\frac{\tilde d_{k+2}/\tilde b_{k+2}-\tilde d_{k+1}/\tilde b_{k+1}}{\tilde d_{k+1}/\tilde b_{k+1}-\tilde d_{k}/\tilde b_{k}}=1.$$
We thus come to the conclusion that condition (B2)$_{ a}$ is satisfied.

{\it Case 2: Suppose that $a,b,\theta$ are not all equal.} In this case, since $(b-a)(b-\theta)\ge0,$ then for some number $k_0>0,$ $\tilde a_k/\tilde b_k=\frac{a+r_k}{b+r_k}+\frac{\theta+r_{k+1}}{b+r_{k+1}},k
\ge k_0$ is monotone in $k.$ That is, $\tilde a_k/\tilde b_k\ne \tilde a_{k+1}/\tilde b_{k+1},k\ge k_0.$

At the first place, if $b+d-a-\theta=0$ and $a\theta-bd=0,$ then it follows from \eqref{fdb} that $\tilde d_k/\tilde b_k=\tilde d_{k+1}/\tilde b_{k+1}=0,\forall k\ge1.$ Using  \eqref{asx}, we get \begin{align}
  \label{aar}\frac{\tilde a_{k+2}}{\tilde b_{k+2}}- \frac{\tilde a_{k+1}}{\tilde b_{k+1}}=\frac{b-a}{b^2}(r_{k+2}-r_{k+1})+\frac{b-\theta}{b^2}(r_{k+3}-r_{k+2})+o(r_{k}^2).
\end{align}
Since $a,b,\theta$ are not all equal and $(b-a)(b-\theta)\ge0,$ we have $(b-a)/b^2+(b-\theta)/b^2\ne0.$
Consequently, we have  $\lim_{k\rto}\frac{\tilde a_{k+2}/\tilde b_{k+2}-\tilde a_{k+1}/\tilde b_{k+1}}{\tilde a_{k+1}/\tilde b_{k+1}-\tilde a_{k}/\tilde b_{k+1}}=1.$ We thus showed that condition (B2)$_{ b}$ holds.

At the second place, suppose that $b+d-a-\theta$ and $a\theta-bd$ are not $0$ simultaneously.
Then from \eqref{fdb} and \eqref{gd}, we have that $\tilde d_k/\tilde b_k,k\ge k_0$ is monotone in $k$ for some number $k_0>0,$ that is,  $\tilde d_k/\tilde b_k\ne \tilde d_{k+1}/\tilde b_{k+1},k\ge k_0.$ But it follows from \eqref{fdb} and \eqref{asx} that
\begin{align*}
  \frac{\tilde d_{k+1}}{\tilde b_{k+1}}-\frac{\tilde d_{k}}{\tilde b_k}&=\z(\frac{b-d}{b^2}-\frac{a(b-\theta)+\theta(b-a)}{b^3}\y)(r_{k+1}-r_k)+o(r_k^2).
\end{align*}
Since $(b-\theta)(b-a)\ge0$ and $a,b,\theta$ are not all equal, taking \eqref{aar} into account, we get
\begin{align*}
&\lim_{k\rto}\frac{\tilde a_{k+2}/\tilde b_{k+2}-\tilde a_{k+1}/\tilde b_{k+1}}{\tilde a_{k+1}/\tilde b_{k+1}-\tilde a_{k}/\tilde b_{k+1}}=1,\\
&\lim_{k\rto}\frac{\tilde d_{k+1}/\tilde b_{k+1}-\tilde d_{k}/\tilde b_k}{\tilde a_{k+1}/\tilde b_{k+1}-\tilde a_{k}/\tilde b_k}=\frac{b(b+d-a-\theta)+2(a\theta-bd)}{b(b-a)+b(b-\theta)}=:\tau
\end{align*}
which is finite and by assumption $\tau \ne \frac{-(a+\theta)
\pm \sqrt{(a+\theta)^2+4(bd-a\theta)}}{2b}.$  Thus we conclude that condition (B2)$_{c}$ holds. Lemma \ref{egc} is proved. \qed
\begin{remark}
  Note that Lemma \ref{egc} excludes the case $a=b=\theta=d,$ which is trivial. Indeed, in this case,
  $M_k=(a+r_k)\left( \begin{array}{cc}
  1& 1 \\
  1 &1 \\
 \end{array}\right)$ and thus $M_1\cdots M_n=(a+r_1)\cdots (a+r_n)\left( \begin{array}{cc}
  2^{n-1}& 2^{n-1} \\
  2^{n-1} &2^{n-1}\\
 \end{array}\right),n\ge1.$

\end{remark}

\vspace{.5cm}

\noindent{{\bf \Large Acknowledgements:}} The authors would like to thank Prof.  Hong, W.M. for introducing to us the basics of BPVE  and  Prof. Vatutin, V. for some discussions on the distribution of the extinction time when writing the paper. Finally the authors were in debt to two
referees who read the paper carefully and gave very good suggestions which help to improve the paper to a large extent. This project is supported by National
Natural Science Foundation of China (Grant No. 11501008).



\begin{thebibliography}{99}
\addtolength{\itemsep}{-0.5em}


\bibitem{bp} \textsc{Bhattacharya, N. and Perlman, M.} (2017).  Time inhomogeneous branching processes conditioned on
non-extinction. {\it arXiv:} 1703.00337.

\bibitem{bcn} \textsc{Biggins, J. D., Cohn, H. and Nerman, O.} (1999). Multi-type branching in varying environment. {\it Stochastic Process. Appl.} \textbf{83}, 357-400.

        \bibitem{cw} \textsc{Cohn, H. and Wang, Q.} (2003). Multitype branching limit behavior. {\it Ann. Appl. Probab.} \textbf{13}, 490-500.


\bibitem{dhkp} \textsc{Dolgopyat, D., Hebbar, P., Koralov, L. and Perlman, M.} (2018). Multi-type branching processes with time-dependent branching rates. {\it J. Appl. Probab.}  \textbf{55},  701-727.


\bibitem{dy99} {\sc Dyakonova, E. E.} (1999). Asymptotic behaviour of the probability of
non-extinction for a multi-type branching
process in a random environment. {\it Discrete Math. Appl.} \textbf{9}, 119-136.

\bibitem{fuj} {\sc Fujimagari, T.} (1980).   On the extinction time distribution of a branching process in
varying environments. {\it Adv. in Appl. Probab.}  \textbf{12},  350-366.
\bibitem{j} {\sc Jagers, P.} (1974).  Galton-Watson processes in varying environment. {\it J. Appl. Probab.} \textbf{11}, 174-178.

    \bibitem{jon}  {\sc Jones, O. D.} (1997). On the convergence of multitype branching processes with varying environments. {\it Ann. Appl. Probab.} \textbf{7}, 772-801.

        \bibitem{ker} {\sc Kersting, G.} (2020). A unifying approach to branching processes in a varying environment. {\it J. Appl. Probab.} \textbf{57}, 196-220.

\bibitem{va} {\sc Kersting, G. and Vatutin, V.} (2017). {\it Discrete time branching processes in random environment.} John Wiley \& Sons, Inc., USA.

        \bibitem{lin} {\sc Lindvall, T.} (1974). Almost sure convergence of branching processes in varying and random environments.
{\it Ann. Probab.} \textbf{2}, 344-346.
        \bibitem{lor} {\sc Lorentzen, L.} (1995). Computation of limit periodic continued fractions. A survey. {\it Numer. Algorithms} \textbf{10}, 69-111.
            \bibitem{lw}  {\sc Lorentzen, L. and Waadeland, H.} (2008). {\it Continued fractions. 2nd. Ed., Volume 1: convergence theory.} Atlantis Press, Paris.

               \bibitem{ms} {\sc Macphee, I. M. and  Schuh, H. J.} (1983). A Galton-Watson branching process in varying
environments with essentially constant means and two rates of growth. {\it Aust. N. Z. J. Stat.} \textbf{25}, 329-338.


     \bibitem{wh} {\sc Sun, H. Y. and Wang, H. M.} (2020). On a maximum of nearest-neighbor random walk with asymptotically zero drift on lattice of positive half line. {\it arXiv:} 2004.12422.
         \bibitem{hs20} {\sc Sun, H. Y. and Wang, H. M.} (2020). Asymptotics of product of nonnegative 2-by-2 matrices  with applications to random walks with asymptotically zero drifts. {\it arXiv:} 2004.13440.


 \end{thebibliography}
\end{document}